\documentclass[11pt]{article}
\usepackage{epsfig}
\usepackage{amssymb}
\usepackage{amsmath}
\usepackage{mathrsfs}
\newcommand{\be}{\begin{equation}}
\newcommand{\ee}{\end{equation}}
\newcommand{\bea}{\begin{eqnarray}}
\newcommand{\eea}{\end{eqnarray}}

\usepackage{bm}
\usepackage{color}
\usepackage{graphicx}
\usepackage{cite}

\begin{document}


\title{On the Sylvester program and Cayley algorithm \\ for vector partition reduction}
\author{Boris Y. Rubinstein\\
Stowers Institute for Medical Research
\\1000 50th St., Kansas City, MO 64110, U.S.A.}
\date{\today}

\maketitle
\begin{abstract}
A vector partition problem asks for a number of nonnegative
integer solutions to a system of linear Diophantine
equations with integer nonnegative coefficients. 
J.J. Sylvester suggested to use variable elimination for the reduction of vector partition to a 
sum of scalar partitions. 
In the simplest case of two equations with positive coefficients A. Cayley performed a reduction 
of the corresponding double partition to a sum of scalar partitions using the algorithm subject to 
a set of conditions on the coefficients. 
We suggested a modification of the 
original Cayley algorithm for the cases when these conditions are not satisfied  \cite{RubDouble2023}.
In this manuscript we generalize the modified Cayley algorithm
to arbitrary number of the Diophantine equations to accomplish the Sylvester program
of the vector partition reduction to a combination of scalar partitions.
\end{abstract}



\section{Introduction}
\label{intro0}

The problem of integer partition into a set of integers 
is equivalent to counting number
of nonnegative integer solutions of the Diophantine equation
\be
s = \sum_{i=1}^m x_i d_i = {\bf x}\cdot{\bf d}.
\label{coin1}
\ee
A scalar partition function $W(s,{\bf d}) \equiv W(s,\{d_1,d_2,\ldots,d_m\})$
solving the above problem is a
number of partitions of an integer $s$ into positive integers
$\{d_1,d_2,\ldots,d_m\}$. The generating function for $W(s,{\bf d})$
has a form
\be
G(t,{\bf d})=\prod_{i=1}^m\frac{1}{1-t^{d_{i}}}
 =\sum_{s=0}^{\infty} W(s,{\bf d})\;t^s\;.
\label{WGF}
\ee
Introducing a notation $C[t^s](f(t))$ for a coefficient of $t^s$ in the 
expansion of a function $f(t)$ we have
\be
W(s,{\bf d}) =
C[t^s]\left( \prod_{i=1}^m(1-t^{d_{i}})^{-1} \right).
\label{const1}
\ee
The first attempt to consider the analytical properties of the function $W(s,{\bf d})$
was made by Cayley in his pioneering work \cite{Cayley1855} where he considered the 
partial fraction decomposition of the generating function $G(t,{\bf d})$.
Cayley found that the scalar partition can be written as a 
sum of a polynomial and expressions of the other type containing the 
factors with integer periods. Based on these findings
Sylvester proved \cite{Sylv2} a theorem about splitting of the 
scalar partition into periodic and non-periodic parts and showed that the
partition function may be presented as a sum of "waves"
\be
W(s,{\bf d}) = \sum_{j=1} W_j(s,{\bf d})\;,
\label{SylvWavesExpand}
\ee
where summation runs over all distinct factors
of the elements of the generator vector ${\bf d}$.
The wave $W_j(s,{\bf d})$ is a quasipolynomial in $s$
closely related to prime roots $\rho_j$ of unity being a coefficient of
${t}^{-1}$ in the series expansion in ascending powers of $t$ of
a function
\be
F_j(s,t)=
\sum_{\rho_j} \frac{\rho_j^{-s} e^{st}}{\prod_{k=1}^{m}
\left(1-\rho_j^{d_k} e^{-d_k t}\right)}\;.
\label{generatorWj}
\ee
The summation is made over all prime roots of unity
$\rho_j=\exp(2\pi i n/j)$ for $n$ relatively prime to $j$
(including unity) and smaller than $j$.
It is possible to
express any Sylvester wave $W_j(s,{\bf d})$ and the scalar partition $W(s,{\bf d})$ 
as a finite sum of the Bernoulli polynomials of
higher order \cite{Rub04}.

Consider a {\it vector partition function} $W({\bf s},{\bf D})$ counting the number of integer
nonnegative
solutions ${\bf x} \ge 0$
to a linear system $ {\bf s} = {\bf D} \cdot {\bf x}$, where
${\bf D}$ is a nonnegative integer $l \times m$ generator matrix ($l < m$).
The vector partition 
is a natural generalization of scalar partition defined in 
(\ref{WGF}-\ref{SylvWavesExpand}) to the vector argument.
The generating function for the vector partition reads
\be
G({\bf t},{\bf D})=\prod_{i=1}^m \frac{1}{1-{\bf t}^{{\bf c}_i}} =
\sum_{{\bf s}} W({\bf s},{\bf D}) {\bf t^s},
\quad
{\bf t^s} = \prod_{k=1}^l t_k^{s_k},
\quad
{\bf t}^{{\bf c}_i} = \prod_{k=1}^l t_k^{c_{ik}},
\label{WvectGF}
\ee
where ${\bf c}_i =\{c_{i1},c_{i2},\ldots,c_{il}\},\ (1 \le i \le m)$
denotes the $i$-th column of the matrix ${\bf D}= \{{\bf c}_1,{\bf c}_2,\ldots,{\bf c}_m\}$.
Generalizing the coefficient notation (\ref{const1})
to the case of function of several variables
find for $W({\bf s},{\bf D})$
\be
W({\bf s},{\bf D}) 
= C[{\bf t}^{\bf s}]
\left(
\prod_{i=1}^m \frac{1}{1-{\bf t}^{{\bf c}_i}}
\right).
\label{GFvect0}
\ee

Several approaches were suggested for vector partition computation including 
the method of residues \cite{Beck2004,Szenes2003} and geometric decomposition
into so called {\em chambers} \cite{Sturmfels1995} -- 
regions in $l$-dimensional space where the partition function
is characterized by a specific expression.
Fundamental results 
on the representation of the generating function (\ref{WvectGF}) into a signed sum of 
simpler expressions of the same type are known \cite{Barvinok1993, Brion1997, DeLoera2004, Milev2023} leading to  
development of the software packages {\it Barvinok} ({\tt https://doc.sagemath.org/html/en/reference/spkg/barvinok.html}) 
and {\it LattE} 
({\tt http://www.math.ucdavis.edu/$\sim$latte/}, page not found as of May 2025) introduced in 1990s and 2000s respectively. 
These programs as well as the new package ``calculator'' ({\tt http://github.com/tmilev/calculator}, 2023) allow effective 
computation of chambers but do not produce explicit expressions for vector partitions in these chambers.
Vector partitions appear to be a powerful tool for the computation of the 
Kronecker coefficients arising as the structure constants in the decomposition
of a tensor product of irreducible representations of the symmetric group
into irreducible representations \cite{Mishna2021},
examples of computation of the Kronecker coefficients using the {\it Barvinok} software
can be found in \cite{Mishna2024}.

An alternative method of vector partition computation was suggested being
a direct generalization of the approach developed in \cite{Rub04}
where {\em vector} Bernoulli and Eulerian polynomials
of higher order were introduced to find explicit
expression for $W({\bf s},{\bf D})$  \cite{Rub06}. A drawback of this
approach is that it does not provide any mechanism to define
the chamber boundaries.
As scalar partition computation requires just known functions (the Bernoulli polynomials of higher order) \cite{Rub04}
it is promising to obtain a reduction method expressing a vector partition
through a combination of scalar ones.
 
The first step in this direction was made in the 19th century by J.J. Sylvester
\cite{Sylv1,Sylv2} who suggested an iterative procedure of 
reduction of a vector partition into a sum of scalar partitions; in \cite{Sylv1} he gave
an example of such a reduction for a double partition. 
The method is based on the elimination one-by-one of the variables $x_i$ 
from the linear system ${\bf s}={\bf D}\cdot{\bf x}$.
Based on this approach A. Cayley 
developed  an algorithm of double partition reduction
and established a set of conditions of the method applicability \cite{Cayley1860}.
Cayley showed that each column ${\bf c}_i$ of the two-row matrix ${\bf D}$ gives rise 
to a scalar partition and elements of this column must be relatively prime (we call it {\it P-column}). 
Another limitation of the Cayley method is that the columns ${\bf c}_i$
should represent {\it noncollinear} vectors.
An application of the Cayley algorithm to the computation of the  Gaussian polynomial
coefficients is considered in \cite{RubSylvCayl}.
A short historical excurse into the (seemingly forgotten) ideas of Sylvester and Cayley
is given in Section \ref{Sylvester-Cayley} while 
the essence of the original Cayley algorithm is presented in Section \ref{Cayleyorig}.

The author of this manuscript considered 
the case when the Cayley formula for double partition fails and discussed 
its modification that successfully resolves the problem \cite{RubDouble2023}.
Specifically, the elimination of a column with the greatest common divisor (GCD)
of its elements $\mbox{gcd}({\bf c}_i)>1$ ({\it NP-column})
leads not to a single scalar partition but to a weighted sum of 
such partitions with modified argument.
This result is described in Section \ref{Cayleynew} and discussed in Section \ref{GCD}.

In Section \ref{Multi} we extend the reduction algorithm 
to $l$-tuple (for arbitrary $l \ge 2$) vector partition defined by a  
$l \times m$ nonnegative generator matrix ${\bf D}$ made of $m$ columns ${\bf c}_i$.
The first step of the algorithm discussed in Section \ref{MultiPFE}
is the partial fraction expansion of 
the generating function (\ref{WvectGF}) into $m$ terms corresponding to 
elimination of each of $m$ unknowns.
The contribution of a column is 
computed in Section \ref{MultiEval} and it is either a 
single partition for P-columns (Section \ref{Multicompact}) 
or a superposition of such partitions (Section \ref{Multiext}) for NP-column. 
When the last element of a column is zero ({\it LZ-column}) both the 
original Cayley method and its generalization cannot be applied directly so that this 
case requires a special treatment. In Section \ref{ZeroMulti} we prove that 
the contribution of such a column to the vector partition is {\it zero}.

The case of collinear columns ({\it C-columns})
of the generator matrix ${\bf D}$ is first addressed in Section \ref{Convolution} 
where we show that 
the vector partition reduces to a convolution of 
scalar partitions. 
An alternative approach of dealing with the C- and NP-columns
was suggested by Sylvester \cite{SylvLectures1859}, namely,
to enlarge the system ${\bf s} = {\bf D} \cdot {\bf x}$ to 
eliminate these columns and then use the variable elimination procedure of partition reduction.
The application of this method is discussed in Section \ref{Extension}.
Additional details of the iterative reduction procedure are discussed 
in Section \ref{Muirhead}    
while the concluding remarks are 
presented in Section \ref{Discuss}.

\section{Sylvester-Cayley method}
\label{Sylvester-Cayley}
The problem of vector (compound) partitions has a long history that started in 1768 
when Euler announced the generating function (\ref{WvectGF}) for double partition ($l=2$)
\cite{Euler1849} and followed by  J.J. Sylvester 
presenting an alternative approach in 1858. Sylvester wrote:
{\it "Any given system of simultaneous simple equations to be solved in positive 
integers being proposed, the determination of the number of solutions of which 
they admit may in all cases be made to depend upon the like determination for 
one or more systems of equations of a certain fixed standard form. When 
a system of $r$ equations between $n$ variables of the aforesaid standard form 
is given, the determination of the number of solutions in positive integers of 
which it admits may be made to depend on the like determination for 
$$\frac{n(n-1)\ldots (n-r+2)}{1\cdot 2 \ldots (r-1)}$$
single {\bf independent} equations derived from those of the given system by the 
ordinary process of elimination, with a slight modification; the result 
being obtained by taking the sum of certain numerical multiples (some positive, 
others negative) of the numbers corresponding to those independent  
determinations. This process admits of being applied in a variety of modes, 
the resulting sum of course remaining unaltered in value whichever mode is employed, only 
appearing for each such mode made up of a different set of component parts"}  \cite{Sylv1}.

Then he added in the footnote:
{\it "If there be $r$ simultaneous simple 
equations between $n$ variables (in which the coefficients are all positive or negative integers) forming 
a definite system (that is, one in which no variable can become indefinitely great in the positive 
direction without one or more of the others becoming negative), and if the $r$ coefficients belonging to 
each of the same variable are exempt from a factor common to them all, and if not more than $r - 1$ 
of the variables can be eliminated simultaneously between the $r$ equations, then the determination 
of the number of positive integer solutions of the given system may be made to depend on like 
determinations for each of $n$ derived independent systems, in each of which the number of variables 
and equations is one less than in the original system"} \cite{Sylv1}.

In other words, Sylvester claimed that $r$-tuple vector partition for $n$ variables can be reduced to 
a sum of $\binom{n}{r-1}$ scalar 
partitions 
and the reduction algorithm is an iterative process based on the 
variable elimination.
Sylvester also described the general conditions of the 
method applicability -- each column should be the coprime set, and 
the elimination of any variable should not lead to the elimination 
of any other variable, i.e., the columns should be noncollinear.
Sylvester considered a specific double partition problem as 
an illustration of his method and determined regions
({\it chambers}) on a plane 
each having a unique
expression for vector partition valid in this region only.
He showed that the 
expressions in the adjacent chambers coincide at their
common boundary (see also \cite{Sturmfels1995}).

Two years later this approach was successfully applied by 
Cayley \cite{Cayley1860} to double partitions subject
to the Sylvester conditions on the matrix ${\bf D}$ elements.
with the added requirement that the inequalities $c_{i2} < s_2+2, \ 1\le i \le m,$ should hold.
In the opening paragraph 
Cayley noticed:
{\it "The subject (as I am aware) has hardly been 
considered except by Professor Sylvester, and it is 
greatly been regretted that only an outline
of his valuable researches has been published: 
the present paper contains the demonstration of a theorem, due to him,
by which (subject to certain restrictions)
the question of Double Partitions is made to depend upon
the ordinary theory of Single Partitions"} \cite{Cayley1860}. 
The manuscript demonstrated that the procedure leads to 
a sum of $m$ scalar partition terms each corresponding to 
elimination of single variable $x_i$ from the vector ${\bf x}$. 

It should be noted that both Sylvester and Cayley usually considered 
coefficients of the system of Diophantine equations ${\bf s} = {\bf D} \cdot {\bf x}$
being positive 
and only 
once Cayley mentioned that the elements of the matrix ${\bf D}$ as well as $s_i$ {\it ``being all positive integer numbers,
not excluding zero''} \cite{Cayley1860} but never explicitly discussed such a case.

Both manuscripts \cite{Sylv1} and  \cite{Cayley1860} were mentioned by L.E. Dickson in the 
then comprehensive volume \cite{DicksonV2} on the theory of numbers
devoted to Diophantine analysis, but 
the literature search shows 
that both the reduction idea of 
Sylvester and the result of Cayley were forgotten.
There are no publications discussing Sylvester research program on 
vector partition reduction to scalar ones. 
In my opinion the ideas and results of Sylvester and Cayley deserve not to be neglected, and this manuscript is an attempt to 
revive them and pave way to a general algorithm for vector partition expression through scalar partitions only.

\section{Cayley method for double partitions}
\label{Cayley}


Computation of vector partition $W({\bf s},{\bf D})$ in (\ref{GFvect0}) is performed by iterative elimination
of $(l-1)$ variables $t_k, (2 \le k \le l)$ \cite {Beck2004}.
Each elimination step includes partial fraction expansion (PFE)
w.r.t. the eliminated variable
with subsequent coefficient evaluation.
This step is equivalent to 
elimination of $k$-th row of augmented $l\times(m+1)$ matrix
${\bf E}= \{{\bf s}={\bf c}_0,{\bf c}_1,{\bf c}_2,\ldots,{\bf c}_m\}$ 
made of the matrix ${\bf D}$
and vector ${\bf s}$; 
the same time the numerical columns of ${\bf D}$
are eliminated one by one. The number of
new $(l-1)\times m$ matrices is equal to $m$.
This algorithm was employed in
\cite{Cayley1860} for a two-row positive matrix 
(see also \cite{Sylv1}). 
The original Cayley algorithm and its generalization 
are exhaustively discussed in \cite{RubDouble2023} and below we 
present their main points.





\subsection{Original Cayley algorithm}
\label{Cayleyorig}

Consider the simplest vector partition case $l=2$ following 
algorithm described in \cite{Cayley1860}.
Denote elements of matrix ${\bf E}$ columns as
${\bf s}={\bf c}_0=\{r,\rho\}^T$ and ${\bf c}_i=\{b_i,\beta_i\}^T,\ (1 \le i \le m),$
where $\{\cdot\}^T$ designates transposition of a vector.
Cayley specified \cite{Cayley1860} the following conditions that
should be met in order to apply the algorithm. First, all fractions 
$b_i/\beta_i$ should be unequal, in other words, the 
columns ${\bf c}_i$ must represent {\it noncollinear} vectors. 
Second, the elements of each column ${\bf c}_i$ should be 
{\it relatively prime}, {\it i.e.}, $\gcd(b_i,\beta_i)=1$.
Finally, all elements
of the second row should satisfy a condition: $\beta_i < \rho+2$.

Assuming all $b_i, \beta_i > 0$
perform PFE step to present
$G({\bf t},{\bf D})$ as sum of $m$ fractions (${\bf t} = \{x,y\}$)
\be
G({\bf t},{\bf D}) =
\prod_{i=1}^m (1-x^{b_i}y^{\beta_i})^{-1} =
\sum_{i=1}^m T_i(x,y),
\quad
T_i = \frac{A_i(x,y)}{1-x^{b_i}y^{\beta_i}},
\label{Cayley01}
\ee
where the functions $A_i$ are rational in $x$ and rational and integral (of degree $\beta_i-1$) in $y$.
The value of $A_1(x,y)$ at $y=y_0=x^{-b_1/\beta_1}$ reads
\be
A_1(x,y_0) =\prod_{i \ne 1}^m (1-x^{b_i}y_0^{\beta_i})^{-1} =
 \prod_{i \ne 1}^m (1-x^{b_i-b_1 \beta_i/\beta_1})^{-1}.
\label{mCayley03}
\ee
Use the relation $1/(1-x) = 1/(1-x^n) \sum_{k=0}^{n-1} x^k$ to write $A_1(x,y_0) $ as
\be
A_1(x,y_0) = \frac{S_1(x,y_0)}{\Pi_1(x)},
\quad
S_1(x,y_0) 
= \sum_{k = 0}^{\beta_1-1} A_{1k}(x) y_0^k,
\quad
\Pi_1(x) =  \prod_{i \ne 1}^m   (1-x^{\beta_1 b_i-b_1 \beta_i}),
\label{mCayley04}
\ee
where $A_{1k}$ are rational functions in $x$.
Consider evaluation of the contribution $C_{r,\rho}^1$ of the term $T_1(x,y)$, it is equal to the coefficient
$$
C_{r,\rho}^1 = C[x^r y^{\rho}](T_1(x,y)) = C[x^r y^{\rho}]
\left(
\frac{A_1(x,y)}{1-x^{b_1}y^{\beta_1}}
\right).
$$
Cayley shows \cite{Cayley1860} that this term evaluates to
\bea
C_{r,\rho}^1 &=&
C[x^r y^{\rho}] \left( \frac{A_1(x,x^{-b_1/\beta_1})}{1-x^{b_1}y^{\beta_1}} \right) 
=
C[x^r] \left(x^{\rho b_1/\beta_1}  A_1(x,x^{-b_1/\beta_1})\right) 
\nonumber  \\
&=&
C[x^{r-\rho b_1/\beta_1}] \left(A_1(x,x^{-b_1/\beta_1})\right)
=C[x^{r\beta_1-\rho b_1}] \left(A_1(x^{\beta_1},x^{-b_1}) \right).
\label{mCayley09a}
\eea
Now from (\ref{mCayley03}) we obtain 
$$
A_1(x^{\beta_1},x^{-b_1}) = \prod_{i \ne 1}^m (1-x^{b_i \beta_1-b_1 \beta_i})^{-1},
$$
and arrive at 
\be
C_{r,\rho}^1 = C[x^{r\beta_1-\rho b_1}] 
\left( \prod_{i \ne 1}^m (1-x^{b_i \beta_1-b_1 \beta_i})^{-1} \right)
=  C[x^{r\beta_1-\rho b_1}] \left(1/\Pi_1(x) \right).
\label{mCayley09}
\ee
Repeating computation for each $T_i(x,y),\ 1 \le i \le m,$ we obtain 
for double partition 
$$
W(\bm s, {\bf D}) = \sum_{i=1}^m C_{r,\rho}^i = 
\sum_{i=1}^m 
C[x^{r\beta_i-\rho b_i}] 
\left( \prod_{j \ne i}^m (1-x^{b_j \beta_i-b_i \beta_j})^{-1} \right).
$$
Recall the definition of the scalar partition function (\ref{const1}) and rewrite it as 
$$
W(s, {\bf d}) = C[t^s] \left( \prod_{i = 1}^m (1-t^{d_i})^{-1} \right),
$$
to obtain an expression of the double partition as a sum of 
scalar partitions
\be
W({\bf s},{\bf D}) = \sum_{i=1}^m W^2_i({\bf s}) =
\sum_{i=1}^m W(L_i({\bf s}),{\bf d}_i),
\quad
L_i({\bf s}) = r\beta_i - b_i\rho,
\quad
d_{ij} = b_j\beta_i-b_i\beta_j,\ j\ne i.
\label{mCayley10}
\ee
This compact expression is the main result of the original Cayley algorithm presented in \cite{Cayley1860}.
Introducing $2 \times 2$ matrices made of the columns of augmented matrix ${\bf E}$
\be
{\bf D}_{i0} = \{{\bf c}_0={\bf s},{\bf c}_i\},
\quad
{\bf D}_{ij} = \{{\bf c}_j,{\bf c}_i\},
\quad
j \ne i,
\label{mCayley010}
\ee
we observe that 
$L_i({\bf s})$ and elements $d_{ij}$ of ${\bf d}_i$ are given by 
determinants of ${\bf D}_{i0}$ and  ${\bf D}_{ij}$, respectively.
Noncollinearity of the columns ${\bf c}_i$ implies that all elements
in the generator sets ${\bf d}_i$ are nonzero, but some of them might be negative
(say, $d_{ij_k}<0$ for $1 \le k \le K_i$). Noting that 
$(1-t^{-a})^{-1} = -t^{a}(1-t^{a})^{-1}$ we find an equivalent scalar partition with positive generators only
\be
W(L_i({\bf s}),{\bf d}_i) =
(-1)^{K} W(L_i({\bf s}) + \sum_{k=1}^{K_i} d_{ij_k},|{\bf d}_i|),
\quad
|{\bf d}_i| = \{|d_{ij}|\}.
\label{mCayley10b}
\ee
The 
partition $W(L_i({\bf s}) + \sum_{k=1}^{K_i} d_{ij_k},|{\bf d}_i|)$ 
in the r.h.s. of (\ref{mCayley10b}) corresponds to the 
Diophantine equation 
$$
\sum_{j\ne i}^m x_j |d_{ij}| = L'_i({\bf s}) = s_1 \beta_i - s_2 b_i +  \sum_{k=1}^{K_i} d_{ij_k},
$$
where $K_i$ integers $d_{ij_k} < 0$ are negative and $x_j \ge 0$ are nonnegative integers.
When $L'_i({\bf s})$ is negative this equation has no solutions and we 
can rewrite (\ref{mCayley10b}) as 
\be
W^2_i({\bf s}) = W(L_i({\bf s}),{\bf d}_i) =
A_i W(L'_i({\bf s}),|{\bf d}_i|),
\quad
A_i = (-1)^{K_i} H(L'_i({\bf s})),
\label{mCayley10b_UnitStep}
\ee
where $H(x)$ is the Heaviside step function defined as 
$$
H(x) = 
\left \{ 
\begin{array}{cc}
0, & x < 0, \\
1, & x \ge 0,
\end{array}
\right .
$$
and this factor determines the chamber to which the scalar partition $W^2_i({\bf s})$ contributes.



The solution (\ref{mCayley10}) is not unique as we can
apply the Cayley method to the matrix ${\bf E}$ with the reverse order of the rows.
It corresponds to the 
elimination of the first row of ${\bf E}$ and gives 
\be
W({\bf s},{\bf D}) = 
\sum_{i=1}^m W(\tilde L_i({\bf s}),\tilde{\bf d}_i),
\quad
\tilde L_i({\bf s}) = -L_i({\bf s}),
\quad
\tilde{\bf d} = -{\bf d}.
\label{mCayley10inv}
\ee
It should be noted that $\tilde L'_i({\bf s}) \ne -L'_i({\bf s})$ and 
the individual contributions $W(L_i({\bf s}),{\bf d}_i)$ and $W(\tilde L_i({\bf s}),\tilde{\bf d}_i)$ differ.
This means that the order of the rows of ${\bf E}$ does not influence the
value of the vector partition $W({\bf s},{\bf D})$ but it
does affect individual summands 
$W^2_i({\bf s})$  in (\ref{mCayley10b_UnitStep}).

It is worth to repeat after Dickson \cite{DicksonV2} that 
the Cayley algorithm corresponds to the elimination of 
each variable $x_i$ from the system ${\bf s}={\bf D}\cdot {\bf x}$ producing 
$m$ equations $L_i = {\bf d}_i \cdot \bm x'_i$ where 
$ \bm x'_i = \{x_1,x_2,\ldots, x_{i-1},x_{i+1},,\ldots,x_m\}$ is 
obtained by dropping $x_i$ from ${\bf x}$.


\subsection{Cayley method generalization}
\label{Cayleynew}

When 
the Cayley method fails 
we suggested
an alternative approach  \cite{RubDouble2023} based on explicit computation of $A_1(x,y)$.
Turning to $S_1(x,y_0)$ in (\ref{mCayley04}) write it as
\be
S_1(x,y_0) =  \prod_{i \ne 1}^m  \frac{(1-x^{\beta_1 b_i}y_0^{\beta_1\beta_i})}{(1-x^{b_i}y_0^{\beta_i})}
= \prod_{i \ne 1}^m 
 \left( 
\sum_{k_i=0}^{\beta_1-1} x^{k_i b_i}y_0^{k_i \beta_i}
\right).
\label{mCayley051}
\ee
Introduce $(m-1)$-dimensional vectors
\be
\bm b'_1 = \{b_2,b_3,\ldots,b_m\},\
\bm \beta'_1 = \{ \beta_2, \beta_3,\ldots,\beta_m\},\
\bm K'_p = \{k_2,k_3,\ldots,k_m\},
\label{mCayley33b}
\ee
with $0 \le k_i \le \beta_1-1$.
Expanding the r.h.s. of (\ref{mCayley051}) we obtain 
\be
S_1(x,y_0) = \sum_{p=1}^{n_1} x^{\bm K'_p \cdot \bm b'_1} y_0^{\bm K'_p \cdot \bm \beta'_1},
\quad
n_1=\beta_1^{m-1}.
\label{mCayley33a}
\ee
For each $\bm K'_p$ we have 
$$
x^{\bm K'_p \cdot \bm b'_1} y_0^{\bm K'_p \cdot \bm \beta'_1} = 
x^{\bm K'_p \cdot \bm b'_1 - b_1 (\bm K'_p \cdot \bm \beta'_1)/\beta_1} 
=x^ {\bm K'_p \cdot \bm b'_1 -  b_1 \lfloor (\bm K'_p \cdot \bm \beta'_1)/\beta_1 \rfloor}
 y_0^{(\bm K'_p \cdot \bm \beta'_1) \bmod{\beta_1}}\;,
$$
where $\lfloor a \rfloor$ denotes the greatest integer less than or equal to real number $a$.
We have for $A_1(x,y)$
\be
A_1(x,y) = \frac{1}{\Pi_1(x)} \sum_{k = 0}^{\beta_1-1} A_{1,j_y}(x) y^{j_y}\;,
\quad
j_y = (\bm K'_p \cdot \bm \beta'_1) \bmod{\beta_1},
\label{mCayley06}
\ee
with 
\be
A_{1,j_y}(x) = \sum_{j_x=N_1^{-}}^{N_1^{+}} a_{1,j_x,j_y} x^{j_x},
\quad
N_1^{-}= \min(\bm K'_p \cdot \bm b'_1 -  b_1 \lfloor (\bm K'_p \cdot \bm \beta'_1)/\beta_1 \rfloor),
\quad
N_1^{+} = (\beta_1-1) \vert  \bm b'_1 \vert,
\label{mCayley06bb}
\ee
where we introduce $L_1$-norm of $m$-component vector $\vert \bm a \vert = \sum_{i = 1}^{m} a_i$.
The integer coefficients $a_{1,j_x,j_y}$ are computed from the relation 
\be
\sum_{j_y=0}^{\beta_1-1}\sum_{j_x=N_1^{-}}^{N_1^{+}} a_{1,j_x,j_y} x^{j_x} y^{j_y} = 
\sum_{p=1}^{n_1} x^{\bm K'_p \cdot \bm b'_1 -  b_1 \lfloor (\bm K'_p \cdot \bm \beta'_1)/\beta_1 \rfloor} 
y^{(\bm K'_p \cdot \bm \beta'_1) \bmod{\beta_1}},
\quad
n_1=\beta_1^{m-1}.
\label{mCayley06bbb}
\ee
Using $A_{1,j_y}(x)$ in expression of $T_1$ and evaluating $C_{r,\rho}^1$ we obtain
\bea
&&C^1_{r,\rho} = \bar W^2_1({\bf s}) = 
\sum_{j_x=N_1^{-}}^{N_1^{+}} a_{1,j_x,j_y}
W(r-j_x-(\rho-j_y)b_1/\beta_1,{\bf d}_1),
\label{mCayley13b} \\
&&j_y = \rho \bmod{\beta_1},
\quad
N_1^{+} = (\beta_1-1) \vert  \bm b'_1 \vert,
\quad
d_{1j} = b_j\beta_1-b_1\beta_j,\ j\ne 1.
\nonumber
\eea 
It is shown in \cite{RubDouble2023} that the term $\bar W^2_i$ 
can be written as a weighted sum of $W_i^2$ with shifted argument
\bea
&&\bar W_i^2({\bf s}) 
= \sum_{{\bf j}={\bf N}_i^{-}}^{{\bf N}_i^{+}} a_{i,{\bf j}} \; \delta_{j_y, \rho \bmod{\beta_i}}
W_i^2(({\bf s} - {\bf j})/\beta_i),
\quad
{\bf j} = \{j_x,j_y\}^T,
\label{mCayley10cc} \\
&&
{\bf N}_i^{-} = \{ \min(\bm K'_p \cdot \bm b'_i -  b_i \lfloor (\bm K'_p \cdot \bm \beta'_i)/\beta_i \rfloor),0\}^T,
\quad
{\bf N}_i^{+} = (\beta_i-1)\{\vert  \bm b'_i \vert,1\}^T,
\nonumber 
\eea
where the vectors $\bm b'_i, \bm \beta'_i$
are given by
$$
\bm b'_i = \{b_1,b_2,\ldots,b_{i-1},b_{i+1},\ldots,b_m\},
\quad
\bm \beta'_i = \{\beta_1,\beta_2,\ldots,\beta_{i-1},\beta_{i+1},\ldots,\beta_m\}.
$$
Using (\ref{mCayley10b_UnitStep}) for each summand in (\ref{mCayley10cc}) 
with ${\bf j} = \{j_1,j_2\}^T$
we can rewrite this expression as
\bea
\bar W_i^2({\bf s}) 
= \sum_{{\bf j}={\bf N}_i^{-}}^{{\bf N}_i^{+}}
A_{i,{\bf j}} W_i^2(({\bf s} - {\bf j})/\beta_i), \ \ 
A_{i,{\bf j}} =  (-1)^{K_i} a_{i,{\bf j}} H(L'_i(({\bf s} - {\bf j})/\beta_i)) \delta_{j_2, s_2 \bmod{\beta_i}} .
\label{mCayley10cc_UnitStep} 
\eea

\subsection{Cayley method for P- and NP-columns}
\label{GCD}


The original Cayley reduction method presented in Section \ref{Cayleyorig} applied to double partition with 
the matrix ${\bf D}$ having $m$ columns generates a sum (\ref{mCayley10}) of $m$ scalar partitions.
This approach fails when at least one of the columns, say the first column
${\bf c}_1$, has GCD of its elements larger than unity $\gcd({\bf c}_1) = \gcd(b_1,\beta_1) = g_1 > 1,$
and we write $b_1=g_1 b^{\star},\ \beta_1=g_1\beta^{\star}$.  
This case should be treated as discussed above in Section \ref{Cayleynew} leading to
$C_{r,\rho}^1$ in (\ref{mCayley13b}).

Consider (\ref{mCayley13b}) and note that when $\gcd(b_1,\beta_1) = g_1 = 1$ 
the symbolic argument in scalar partitions in the summand is an integer only if
the condition $j_y = \rho \bmod{\beta_1}$ holds.
This means that when ${\bf c}_1$ is a P-column the condition $j_y = \rho \bmod{\beta_1}$ 
might be dropped and {\it all} terms in the sum in (\ref{mCayley13b}) contribute to $C_{r,\rho}^1$
that has an alternative compact form used in (\ref{mCayley10}).

In case of NP-column ${\bf c}_1$ with $g_1 > 1$ 
we have 
$
(\rho-j_y)b_1/\beta_1 = (\rho-j_y)b^{\star}/\beta^{\star},
$
and for given value of $\rho$ there exist exactly $g_1$ values of $j_y$ 
for which the scalar partition symbolic argument $r-j_x-(\rho-j_y)b_1/\beta_1$
has integer value. 
The condition $j_y = \rho \bmod{\beta_1}$  
selects only one of $g_1$ values of $j_y$ that should be used for the 
estimate of the contribution to the vector partition. 
Thus {\it not all} (but each $g_1$-th only) summands with integer scalar partition symbolic argument 
contribute and this fact prevents the compactification of the superposition  (\ref{mCayley13b})
into a single term. An example of NP-column reduction in double partition is discussed in
Appendix \ref{Triple1}.


Note that for $\beta_1=0$ the 
expression for $S_1(x,y_0)$ in (\ref{mCayley051}) vanishes  implying
vanishing contribution $C^1_{r,\rho}= 0$. 
We also observe that the result (\ref{mCayley10b}) with $i=1$ and $\beta_1=0$ formally
produces all $d_{1j} < 0$ and $L_1({\bf s}) = -s_2 b_1$ 
leading to negative $L'_1({\bf s}) = -s_2 b_1 + \sum_{j=2}^{m} d_{1j}$ for any 
{\it positive} $s_2$. This means that the scalar partition 
$W(L_1({\bf s}),{\bf d}_1) = W(L'_1({\bf s}), |{\bf d}_1|)$ equals zero.
Thus the rule of thumb for double partitions with nonnegative generators $b_i,\beta_i$ and 
$r,\rho \ge 0$ is -- the contribution of the LZ-column with the last zero $\beta_i=0$ 
is zero.

The same time we have to recall that in the course of the derivation of (\ref{mCayley10b})
we explicitly assumed that $\beta_1 \ne 0$. Nevertheless it is possible to show that 
in case $\beta_1=0$ the contribution $C^1_{r,\rho}$ indeed vanishes. 
It appears to be a particular case of the general result valid for a $l$-tuple partition ($l \ge 2$)
the first LZ-column ${\bf c}_1$ of the $l \times m$ matrix
${\bf D}$. The reader is referred to Section \ref{ZeroMulti} where the details of the derivation are presented.


\section{$l$-tuple partition reduction}
\label{Multi}

Consider a partition $W({\bf s},{\bf D})$ for a 
non-negative  $l \times m\ (l < m)$ matrix ${\bf D}=\{{\bf c}_1,{\bf c}_2,\ldots,{\bf c}_m\}$
made of $m$ columns 
${\bf c}_i=\{c_{i1},c_{i2},\ldots,c_{il}\}^T$
in assumption of column noncollinearity.
We would also use the equivalent notation 
$W({\bf E})\equiv W({\bf s},{\bf D})$ for the partition with corresponding augmented matrix 
${\bf E}=\{{\bf c}_0={\bf s},{\bf c}_1,{\bf c}_2,\ldots,{\bf c}_m\}$ 
obtained by adding the symbolic argument column ${\bf s}$ to matrix ${\bf D}$.


\subsection{Partial fraction expansion}
\label{MultiPFE}

Perform a PFE step (elimination of the last $l$-th row of ${\bf D}$) in assumption that all elements of this row are {\it positive} to present
$G({\bf t},{\bf D})$ with ${\bf t} = \{t_{1},t_{2},\ldots,t_{l}\}$ as sum of $m$ fractions
using   
additional shortcut notations 
\be
\hat{\bf t}=\{t_{1},\ldots,t_{l-1}\},
\quad
\hat{\bf c}_i=\{c_{i1},\ldots,c_{i,l-1}\},
\quad
\hat{\bf s}=\{s_{1},\ldots,s_{l-1}\},
\quad
\hat{\bf j}=\{j_{1},\ldots,j_{l-1}\}.
\nonumber
\ee
We have
\be
G({\bf t},{\bf D}) =
\prod_{i=1}^m (1-\hat{\bf t}^{\hat{\bf c}_i}t_l^{c_{il}})^{-1} =
\sum_{i=1}^m T_i(\hat{\bf t},t_l),
\quad
T_i = \frac{A_i(\hat{\bf t},t_{l})}{1-\hat{\bf t}^{\hat{\bf c}_i}t_l^{c_{il}}},
\label{mCayley101}
\ee
where functions $A_i$ are rational in $\hat{\bf t}$ and rational and integral (of degree $c_{il}-1$) in $t_l$.
Consider the term $A_1(\hat{\bf t},t_l)$ that corresponds to elimination of the last row and the 
first column ${\bf c}_1$ of ${\bf E}$. Introduce
matrix ${\bf E}_1$ produced from ${\bf E}$ by dropping the column ${\bf c}_1$
and present it as an array of $l$ rows ${\bf r}_{1k},\ 1 \le k \le l,$ each having $m$ elements
${\bf E}_1 = \{{\bf r}_{11},{\bf r}_{12},\ldots,{\bf r}_{1l}\}$.
The term $A_1(\hat{\bf t},t_l)$ for 
$t_{l0} = \hat{\bf t}^{-\hat{\bf c}_1/c_{1l}}$
turns into 
\be
A_1(\hat{\bf t},t_l) = \prod_{i \ne 1}^m (1-\hat{\bf t}^{\hat{\bf c}_i}t_l^{c_{il}})^{-1},
\label{mCayley102}
\ee
that is 
\be
A_1(\hat{\bf t},\hat{\bf t}^{-\hat{\bf c}_1/c_{1l}}) = 
\prod_{i \ne 1}^m (1-\hat{\bf t}^{\hat{\bf c}_i-c_{il}\hat{\bf c}_1/c_{1l}})^{-1}.
\label{mCayley103}
\ee
Following the approach of Section \ref{Cayleyorig} 
write (\ref{mCayley103}) as
\be
A_1(\hat{\bf t},\hat{\bf t}^{-\hat{\bf c}_1/c_{1l}}) = 
S_1(\hat{\bf t},t_{l0})/\Pi_1(\hat{\bf t}),
\quad
\Pi_1(\hat{\bf t}) = \prod_{i \ne 1}^m  (1-\hat{\bf t}^{c_{1l}\hat{\bf c}_i}t_{l0}^{c_{1l}c_{il}}) =
 \prod_{i \ne 1}^m   (1-\hat{\bf t}^{c_{1l}\hat{\bf c}_i-c_{il}\hat{\bf c}_1}).
\label{mCayley103a}
\ee
The numerator reads
\be
S_1(\hat{\bf t},t_{l0}) =
 \prod_{i \ne 1}^m 
 \left( 
\sum_{k_i=0}^{c_{1l}-1} \hat{\bf t}^{k_i \hat{\bf c}_i}t_{l0}^{k_i c_{il}}
\right)
= \sum_{j_l = 0}^{c_{1l}-1} A_{1,j_l}(\hat{\bf t}) t_{l0}^{j_l},
\label{mCayley3051}
\ee
where $A_{1,j_l}$ are rational functions in $\hat{\bf t}$.
Define $l$ vectors $\bm r'_{1k},\ 1\le k \le l,$ 
by retaining last $(m-1)$ (only numerical) elements of vectors ${\bf r}_{1k}$,
and a set of $c_{1l}^{m-1}$ vectors $\bm K'_p = \{k_2,k_3,\ldots,k_m\}$ 
with $0 \le k_i \le c_{1l}-1$.
Expanding the product in (\ref{mCayley3051}) we obtain a sum
\be
S_1(\hat{\bf t},t_{l0}) = \sum_{p=1}^{n_1} 
t_{l0}^{\bm K'_p \cdot \bm r'_{1l}}
\prod_{k=1}^{l-1}
t_{k}^{\bm K'_p \cdot \bm r'_{1k}} ,
\quad
n_1=c_{1l}^{m-1},
\quad
t_{l0}=\hat{\bf t}^{-\hat{\bf c}_1/c_{1l}}.
\label{mCayley333a}
\ee
For each $\bm K'_p$ we have 
$$
t_{l0}^{\bm K'_p \cdot \bm r'_{1l}} = 
\prod_{k=1}^{l-1}
t_{k}^{-c_{1k}(\bm K'_p \cdot \bm r'_{1k})/c_{1l}}
$$
and obtain
$$
t_{l0}^{\bm K'_p \cdot \bm r'_{1l}}
\prod_{k=1}^{l-1}
t_{k}^{\bm K'_p \cdot \bm r'_{1k}}  =
\prod_{k=1}^{l-1}
t_{k}^{\bm K'_p \cdot \bm r'_{1k} - c_{1k}(\bm K'_p \cdot \bm r'_{1l})/c_{1l}} = 
\prod_{k=1}^{l-1}
t_{k}^{\nu_k},
$$
where the rational exponent $\nu_k$ of $t_k$ is split into an integer $j_k$ and 
a fractional part. The fractional part of $\nu_k$ 
gives an exponent $j_l$ of $t_{l0}$,
while $t_k^{j_k}$ contributes to $A_{1,j_l}(\hat{\bf t})$ in (\ref{mCayley3051})
as a factor in $\hat{\bf t}^{\hat{\bf j}}$.
We observe that 
\be
j_l = (\bm K'_p \cdot \bm  r'_{1l}) \bmod{c_{1l}},
\label{mCayley333d}
\ee
and obtain
\be
\bm K'_p \cdot \bm r'_{1l} = c_{1l} t + j_l,
\
t = \lfloor (\bm K'_p \cdot \bm r'_{1l})/c_{1l} \rfloor,
\
0 \le t=T_1(j_l) \le \lfloor (R_1-j_l)/c_{1l} \rfloor,
\  R_1 = (c_{1l}-1) \vert  \bm r'_{1l} \vert.
\label{mCayley333dd}
\ee 
This leads to  
\be
t_{l0}^{\bm K'_p \cdot \bm r'_{1l}}
\prod_{k=1}^{l-1}
t_{k}^{\bm K'_p \cdot \bm r'_{1k}}  =
t_{l0}^{j_l} \hat{\bf t}^{\hat{\bf j}} = 
t_{l0}^{(\bm K'_p \cdot \bm  r'_{1l}) \bmod{c_{1l}}}
\prod_{k=1}^{l-1}
t_{k}^{j_k},
\label{mCayley333e}
\ee
with
\be
j_k = \bm K'_p \cdot \bm r'_{1k} - c_{1k}\lfloor (\bm K'_p \cdot \bm r'_{1l})/c_{1l} \rfloor = 
\bm K'_p \cdot \bm r'_{1k} -  c_{1k}t,
\quad
1 \le k \le l-1.
\label{mCayley345}
\ee
We arrive at the relation 
\be
\sum_{j_l = 0}^{c_{1l}-1} A_{1,j_l}(\hat{\bf t}) t_{l0}^{j_l} = 
\sum_{p=1}^{n_1} 
t_{l0}^{(\bm K'_p \cdot \bm  r'_{1l}) \bmod{c_{1l}}}
\prod_{k=1}^{l-1}
t_{k}^{\bm K'_p \cdot \bm r'_{1k} - c_{1k}\lfloor (\bm K'_p \cdot \bm r'_{1l})/c_{1l} \rfloor},
\label{mCayley306b}
\ee
that allows to write for $A_{1,j_l}(\hat{\bf t})$
\be
A_{1,j_l}(\hat{\bf t}) = 
\sum_{\hat{\bf j}=\hat{\bf N}_1^{-}}^{\hat{\bf N}_1^{+}} a_{1,\hat{\bf j},j_l} \hat{\bf t}^{\hat{\bf j}},
\quad
\hat{\bf N}_1^{\pm} = \{N_{11}^{\pm},N_{21}^{\pm},\ldots,N_{l-1,1}^{\pm}\}^T,
\label{mCayley306bb}
\ee
with 
$$
N_{k1}^{-} = \min(\bm K'_p \cdot \bm r'_{1k} - c_{1k}\lfloor (\bm K'_p \cdot \bm r'_{1l})/c_{1l} \rfloor),
\quad
N_{k1}^{+} = (c_{1l}-1) \vert  \bm r'_{k1} \vert,
\quad
1 \le k \le l-1.
$$
The expansion coefficients $a_{1,\hat{\bf j},j_l}$ in (\ref{mCayley306bb}) are computed 
directly from 
(\ref{mCayley306b}). 
The relation (\ref{mCayley101}) leads to
\be
T_1(\hat{\bf t},t_l) = \frac{A_1(\hat{\bf t},t_l)}{1-\hat{\bf t}^{\hat{\bf c}_i}t_l^{c_{il}}} = 
\frac{1}{\Pi_1(\hat{\bf t})(1-\hat{\bf t}^{\hat{\bf c}_1}t_l^{c_{1l}})} 
\sum_{j_l = 0}^{c_{1l}-1} A_{1,j_l}(\hat{\bf t}) t_l^{j_l} \;.
\label{mCayley306c}
\ee
Substitute (\ref{mCayley306bb}) in (\ref{mCayley306c}) to arrive at 
\be
T_1({\bf t}) = 
\frac{1}{\Pi_1(\hat{\bf t})(1-{\bf t}^{{\bf c}_1})} 
\sum_{{\bf j}={\bf N}_1^{-}}^{{\bf N}_1^{+}} a_{1,{\bf j}} {\bf t}^{\bf j},
\label{mCayley306gen}
\ee
where 
${\bf N}_1^{\pm}$ are obtained by adding $N_{l1}^{\pm}$ to $\hat{\bf N}_1^{\pm}$
with  $N_{l1}^{-}=0,$ and $N_{l1}^{+}=(c_{1l}-1) \vert  \bm r'_{l1} \vert$.


\subsection{Contribution evaluation}
\label{MultiEval}

Use the definition (\ref{GFvect0}) to compute the contribution 
$C_{\bf s}^1 = C[{\bf t}^{\bf s}](T_1({\bf t}))$ of the term $T_1$
defined in (\ref{mCayley306c})
\be
C_{\bf s}^1 = 
C[\hat{\bf t}^{\hat{\bf s}}]
\left(
\Pi_1^{-1}(\hat{\bf t})
C[t_l^{s_l}]
\left(
\frac{1}{(1-\hat{\bf t}^{\hat{\bf c}_1}t_l^{c_{1l}})}
\sum_{j_l = 0}^{c_{1l}-1} A_{1,j_l}(\hat{\bf t}) t_l^{j_l}
\right)
\right).
\label{Gen401}
\ee
First consider the inner term
\be
C_l^1 = C[t_l^{s_l}]
\left(
\frac{\sum_{j_l = 0}^{c_{1l}-1} A_{1,j_l}(\hat{\bf t}) t_l^{j_l}}{(1-\hat{\bf t}^{\hat{\bf c}_1}t_l^{c_{1l}})}
\right)
= 
\sum_{p=0}^{\infty}
\sum_{j_l = 0}^{c_{1l}-1} 
A_{1,j_l}(\hat{\bf t})\hat{\bf t}^{p\hat{\bf c}_1}
C[t_l^{s_l}]
\left(
 t_l^{pc_{1l}+j_l}
\right),
\label{Gen402}
\ee
and observe that $C[t_l^{s_l}]\left(t_l^{pc_{1l}+j_l}\right)=1$ only 
for $pc_{1l}+j_l = s_l$ otherwise it vanishes. Find $p = (s_l-j_l)/c_{1l}$ and use (\ref{mCayley3051})
to obtain
\be
C_l^1 = \sum_{j_l = 0}^{c_{1l}-1} 
A_{1,j_l}(\hat{\bf t})\hat{\bf t}^{(s_l-j_l)\hat{\bf c}_1/c_{1l}} = 
\hat{\bf t}^{s_l\hat{\bf c}_1/c_{1l}} \sum_{j_l = 0}^{c_{1l}-1} 
A_{1,j_l}(\hat{\bf t}) t_{l0}^{j_l}
=  \hat{\bf t}^{s_l\hat{\bf c}_1/c_{1l}}S_1(\hat{\bf t},t_{l0}).
\label{Gen403}
\ee


\subsubsection{Compact form}
\label{Multicompact}

Substitute (\ref{Gen403}) in (\ref{Gen401}) to find
$$
C_{\bf s}^1 = 
C[\hat{\bf t}^{\hat{\bf s}}]
\left(
\Pi_1^{-1}(\hat{\bf t})\hat{\bf t}^{s_l\hat{\bf c}_1/c_{1l}}S_1(\hat{\bf t},t_{l0})
\right) = 
C[\hat{\bf t}^{\hat{\bf s}-s_l\hat{\bf c}_1/c_{1l}}]
\left(
\Pi_1^{-1}(\hat{\bf t})
S_1(\hat{\bf t},t_{l0}) 
\right)
= 
C[\hat{\bf t}^{\hat{\bf s}-s_l\hat{\bf c}_1/c_{1l}}]
\left(
A_1(\hat{\bf t},t_{l0}) 
\right).
$$
Recall that $t_{l0} =  \hat{\bf t}^{-\hat{\bf c}_1/c_{1l}}$ and 
using (\ref{mCayley103},\ref{mCayley103a}) find 
\be
C_{\bf s}^1 = 
C[\hat{\bf t}^{\hat{\bf s}-s_l\hat{\bf c}_1/c_{1l}}]
\left(
A_1(\hat{\bf t},\hat{\bf t}^{-\hat{\bf c}_1/c_{1l}}) 
\right) = 
C[\hat{\bf t}^{c_{1l}\hat{\bf s}-s_l\hat{\bf c}_1}]
\left(
A_1(\hat{\bf t}^{c_{1l}},\hat{\bf t}^{-\hat{\bf c}_1}) 
\right) 
= C[\hat{\bf t}^{c_{1l}\hat{\bf s}-s_l\hat{\bf c}_1}]
\left(
1/\Pi_1(\hat{\bf t})
\right),
\label{Gen405}
\ee
Employing expression (\ref{GFvect0}) for $l$-tuple partition $W({\bf s},{\bf D})$
we observe 
\be
C_{\bf s}^1 = W^l_1 = W({\bf s}^1,{\bf D}_1),
\ 
{\bf s}^1 = c_{1l}\hat{\bf s}-s_l\hat{\bf c}_1,
\
{\bf D}_1=\{{\bf c}^1_2,{\bf c}^1_3,\ldots,{\bf c}^1_m\},
\
{\bf c}^1_j=c_{1l}\hat{\bf c}_j-c_{jl}\hat{\bf c}_1,\ 2\le j \le m.
\label{Gen407}
\ee
We  observe that the expressions for ${\bf s}^1$ and ${\bf c}^1_j$
follow the same pattern and are just $2 \times 2$ determinants similar to
those introduced in (\ref{mCayley010}). The determinants are computed using the 
columns ${\bf c}^1$ and corresponding columns of the augmented matrix ${\bf E}_1$ 
introduced in the beginning of this Section.


It should be reminded here that the final goal of the Sylvester program 
is the vector partition reduction to a set of scalar partitions and thus the elimination of a single (last)
row of $l \times m$ matrix is just the first step of the procedure. 
We discuss the computational aspects of the Sylvester program further in
Section \ref{Muirhead} but now we have to underline that the reduction procedure 
described above is valid only when the matrix under reduction 
has {\it its last row void of nonpositive elements}.
It follows from (\ref{Gen407}) that some columns ${\bf c}^1_j$
might have negative elements $c^1_{jl} < 0$ in the last position.

A simple procedure similar to one leading 
to (\ref{mCayley10b}) allows to resolve this issue.
Namely, apply the relation 
\be
\left(1-{\bf t}^{{\bf c}}\right)^{-1} = -{\bf t}^{-{\bf c}}\left(1-{\bf t}^{-{\bf c}}\right)^{-1},
\label{mCayley118}
\ee
to all $K_1$  columns ${\bf c}^1_j$ with $c^1_{jl} < 0$ to
update the zeroth column ${\bf s}^1$ and 
the partition itself:
\be
\tilde{\bf c}^1_j = -{\bf c}^1_j,
\quad
\tilde{\bf s}^1 = {\bf s}^1 + \sum_{j=1}^{K_1} {\bf c}^1_j,
\quad
W(\{{\bf s}^1,{\bf c}^1_1,\ldots,{\bf c}^1_m\}) = 
(-1)^{K_1} H(\tilde s^1_{l-1}) W(\{\tilde{\bf s}^1,\tilde{\bf c}^1_1,\ldots,\tilde{\bf c}^1_m\}),
\label{Gen407lastrowpositive}
\ee
where $\tilde{\bf c}^1_i$ for the columns not affected by the 
procedure is equal to its original value $\tilde{\bf c}^1_i = {\bf c}^1_i$.
The factor $H(\tilde s^1_{l-1})$ reflects the fact that the $(l-1)$-th
element of the column $\tilde{\bf s}^1$ should be nonnegative $\tilde s^1_{l-1} \ge 0$.


Note that (\ref{Gen403}) does not work when GCD
of the column ${\bf c}_1$ elements is larger than unity:
$\gcd({\bf c}_1) \equiv \gcd(c_{11},c_{12},\ldots,c_{1l}) = g_1 > 1$. 
Then we have for such 
NP-column ${\bf c}_1=g_1 {\bf c}^{\star}$
leading to $\hat{\bf c}_1/c_{1l} = \hat{\bf c}^{\star}/c_{l}^{\star}$ and 
$t_{l0}^{j_l} =  \hat{\bf t}^{-j_l\hat{\bf c}_1/c_{1l}} = \hat{\bf t}^{-j_l\hat{\bf c}^{\star}/c_{l}^{\star}}$
enters the sum in (\ref{Gen403}) where $1 \le j_l \le c_{1l}-1$, but $ c_{1l}-1 > c_{l}^{\star}$.
The summation does not evaluate into $S_1(\hat{\bf t},t_{l0})$ and thus both (\ref{Gen403}) 
and (\ref{Gen407}) fail.
The case of NP-column requires an approach that leads to 
an expanded form of $C_{\bf s}^1$ contribution.


\subsubsection{General form}
\label{Multiext}

To obtain a general form of the contribution $C_{\bf s}^1$  note that the relation $pc_{1l}+j_l = s_l$
restricts $j_l$ value to $j_l = s_l \bmod{c_{1l}}$ and the sum in (\ref{Gen403})
reduces to a single term
\be
C_l^1 = 
A_{1,j_l}(\hat{\bf t})\hat{\bf t}^{(s_l-j_l)\hat{\bf c}_1/c_{1l}},
\quad
j_l = s_l \bmod{c_{1l}},
\label{Gen408}
\ee
where $A_{1,j_l}(\hat{\bf t})$ is given by (\ref{mCayley306bb})
and we have
\bea
C_{\bf s}^1 &=& 
C[\hat{\bf t}^{\hat{\bf s}}]
\left(
\Pi_1^{-1}(\hat{\bf t})C_l^1
\right) = 
\!\!\sum_{\hat{\bf j}=\hat{\bf N}_1^{-}}^{\hat{\bf N}_1^{+}} \!\! a_{1,\hat{\bf j},j_l}
C[\hat{\bf t}^{\hat{\bf s}}]
\left(
\Pi_1^{-1}(\hat{\bf t}) \hat{\bf t}^{\hat{\bf j}+(s_l-j_l)\hat{\bf c}_1/c_{1l}}
\right)
\nonumber \\
&=&
\!\!\sum_{\hat{\bf j}=\hat{\bf N}_1^{-}}^{\hat{\bf N}_1^{+}} \!\! a_{1,\hat{\bf j},j_l}
C[\hat{\bf t}^{\hat{\bf s}-\hat{\bf j}-(s_l-j_l)\hat{\bf c}_1/c_{1l}}]
\left(
1/\Pi_1(\hat{\bf t})
\right),
\quad
j_l = s_l \bmod{c_{1l}}.
\label{Gen409}
\eea
Consider a summand in (\ref{Gen409}) for some arbitrary value of column $\hat{\bf j}$,
it reduces to
$$
C[\hat{\bf t}^{\hat{\bf s}-\hat{\bf j}-(s_l-j_l)\hat{\bf c}_1/c_{1l}}]
\left(
1/\Pi_1(\hat{\bf t})
\right) = 
W({\bf s}^{1}({\bf j}),{\bf D}_1),
\quad
{\bf s}^{1}({\bf j}) = \hat{\bf s}-\hat{\bf j}-(s_l-j_l)\hat{\bf c}_1/c_{1l},
$$
and write
\be
C_{\bf s}^1 = \bar W^l_1 = 
\!\!\sum_{\hat{\bf j}=\hat{\bf N}_1^{-}}^{\hat{\bf N}_1^{+}} \!\! a_{1,\hat{\bf j},j_l}
W(\hat{\bf s}-\hat{\bf j}-(s_l-j_l)\hat{\bf c}_1/c_{1l},{\bf D}_1),
\quad
j_l = s_l \bmod{c_{1l}}.
\label{Gen410}
\ee
It is convenient to 
rewrite $\bar W^l_1$ using the 
columns ${\bf N}_1^{\pm}$
employed in (\ref{mCayley306gen}) 
\be
\bar W^l_1 = 
\!\!\sum_{{\bf j}={\bf N}_1^{-}}^{{\bf N}_1^{+}} \!\! a_{1,{\bf j}}
W({\bf s}^{1}({\bf j}),{\bf D}_1)
\delta_{j_l,s_l \bmod{c_{1l}}},
\quad
{\bf s}^{1}({\bf j}) = \frac{{\bf s}^{1}-\hat{\bf j}c_{1l}+j_l\hat{\bf c}_1}{c_{1l}}.
\label{Gen411}
\ee

We 
note here that (\ref{Gen411}) in most cases is optional 
but it is compulsory for the NP-column with $g_1 = \gcd({\bf c}_1) > 1$
and its special case of {\it Z-column} with $\hat{\bf c}_1=0,\ g_1=c_{1l}>1$. 


For the NP-column we observe from (\ref{Gen407}) that {\it all} columns 
${\bf c}^1_j$ of ${\bf D}_1$ are proportional to $g_1$. 
This means that the partition 
$W(\hat{\bf s}-\hat{\bf j}-(s_l-j_l)\hat{\bf c}_1/c_{1l},{\bf D}_1)$ in (\ref{Gen410})
should vanish when $g_1$ does not divide its first argument.
It is obvious that with $j_l = s_l \bmod{c_{1l}}$ the last term $(s_l-j_l)\hat{\bf c}_1/c_{1l}$
is proportional to $g_1$ and 
this means that the summands in (\ref{Gen410}) are nonzero only if 
a set of conditions $(\hat{\bf s} - \hat{\bf j}) \bmod{g_1} = 0$  is satisfied.
In order to prevent application of general form (\ref{Gen411}) in further steps of the 
reduction algorithm it is recommended to use its equivalent version 
\be
\bar W^l_1 = 
\!\!\sum_{{\bf j}={\bf N}_1^{-}}^{{\bf N}_1^{+}} \!\! a_{1,{\bf j}}
W({\bf s}^{1}({\bf j})/g_1,{\bf D}_1/g_1)
\delta_{j_l,s_l \bmod{c_{1l}}},
\
(\hat{\bf s} - \hat{\bf j}) \bmod{g_1} = 0,
\
{\bf s}^{1}({\bf j}) = \hat{\bf s}-\hat{\bf j}-(s_l-j_l)\hat{\bf c}_1/c_{1l}.
\label{Gen411GCD}
\ee

Finally performing the removal of $K_1$ negative elements from the last row
of the matrix ${\bf D}_1$ similar to one discussed in Section \ref{Multicompact}
we obtain another form of (\ref{Gen411GCD}) 
\bea
&&\bar W^l_1 = 
\!\!\sum_{{\bf j}={\bf N}_1^{-}}^{{\bf N}_1^{+}} \!\! (-1)^{K_1} a_{1,{\bf j}} H(\tilde s^1_l({\bf j}))
W(\tilde{\bf s}^{1}({\bf j})/g_1,\tilde{\bf D}_1/g_1)
\delta_{j_l,s_l \bmod{c_{1l}}},
\label{Gen411GCDlastrowpositive}\\
&&
\tilde{\bf s}^1({\bf j}) = {\bf s}^1({\bf j}) + \sum_{j=1}^{K_1} {\bf c}^1_j,
\quad
\tilde{\bf D}_1 = \{\tilde{\bf c}^1_1,\ldots,\tilde{\bf c}^1_m\}),
\nonumber
\eea
where $\tilde{\bf c}^1_j$ are defined in 
(\ref{Gen407lastrowpositive}). The vector partition 
$W(\tilde{\bf s}^{1}({\bf j})/g_1,\tilde{\bf D}_1/g_1)$ is ready
for the next reduction step.

For the Z-column we have $\hat{\bf c}_1={\bf 0}$ and 
the argument $\hat{\bf s}-\hat{\bf j}-(s_l-j_l)\hat{\bf c}_1/c_{1l}$ in (\ref{Gen410})
reduces to $\hat{\bf s}-\hat{\bf j}$.
All columns of ${\bf D}_1$ are proportional to $g_1=c_{1l}$
and we deduce that  (\ref{Gen411}) turns into
\be
\bar W^l_1 = 
\!\!\sum_{{\bf j}={\bf N}_1^{-}}^{{\bf N}_1^{+}} \!\! a_{1,{\bf j}}
W((\hat{\bf s}-\hat{\bf j})/c_{1l},{\bf D}_1/c_{1l})
\delta_{j_l,s_l \bmod{c_{1l}}},
\quad
(\hat{\bf s} - \hat{\bf j}) \bmod{c_{1l}} = 0,
\label{Gen411zero}
\ee
that also should be converted to a form similar to (\ref{Gen411GCDlastrowpositive}).


\subsection{Vector partitions with LZ-columns}
\label{ZeroMulti}

Consider evaluation of a vector partition having the first LZ-column with $c_{1l}=0$. 
Rewrite the system ${\bf s} = {\bf D} \cdot {\bf x}$ as 
\be
 {\bf s}  = {\bf c}_1 x_1 + \sum_{i=2}^m {\bf c}_i x_i
\Rightarrow
 {\bf s}'_1 = \sum_{i=2}^m {\bf c}_i x_i,
\quad
{\bf s}'_1 = {\bf s} - {\bf c}_1 x_1 =  \{s_1-c_{11} x_1,s_2-c_{12} x_1,\ldots,s_l\}^T.
\label{Zero1a}
\ee
Introduce $\bm x'_1 = \{x_2,x_3,\ldots,x_m\}$ and use 
the vectors $\bm r'_{1k},\ 1\le k \le l,$ defined in Section \ref{MultiPFE} to obtain 
a system of $l$ equations 
$$
\bm r'_{1k} \cdot \bm x'_1 = s_k-c_{1k} x_1, \ \ (1\le k \le l-1),
\quad
\bm r'_{1l} \cdot \bm x'_1 = s_l.
$$
The matrix $\bm D'_1 =  \{\bm r'_{11},\bm r'_{12},\ldots,\bm r'_{1l}\}$ made of the 
rows $\bm r'_{1k}$ together with the column ${\bf s}'_1$
define a vector partition $W({\bf s}'_1,\bm D'_1)$ for each nonnegative 
value $x_1 \ge 0$.
This means that 
\be
W({\bf s},{\bf D}) = \sum_{x_1=0}^{\infty} W({\bf s}'_1,\bm D'_1), 
\label{Zero13}
\ee
where we assume that the last row of $\bm D'_1$ is free of zeros.
The relations (\ref{Gen407}) determine the 
contribution $W_i^l({\bf s})$ obtained by elimination of the P-column ${\bf c}_i,\ i>1$
as $W({\bf s}^i,{\bf D}_i)$ with ${\bf D}_i=\{{\bf c}^i_1,{\bf c}^i_2,\ldots,{\bf c}^i_m\}$ 
and ${\bf c}^i_j=c_{il}\hat{\bf c}_j-c_{jl}\hat{\bf c}_i,\  1 \le j\ne i \le m,$ 
where $\hat{\bf c}_i$ denotes the column ${\bf c}_i$ without its last element $c_{il}$.
The explicit expression for ${\bf c}^i_1$ reads ${\bf c}^i_1=c_{il}\hat{\bf c}_1-c_{1l}\hat{\bf c}_i = c_{il}\hat{\bf c}_1$.
The column ${\bf s}^i$ is evaluated to ${\bf s}^i = c_{il}\hat{\bf s}-s_l\hat{\bf c}_i$. 

Consider the elimination of the P-column ${\bf c}_1^i,\ i>1$ from the matrix $\bm D'_1$ 
in the partition $W({\bf s}'_1,\bm D'_1)$. It is represented by a single term 
\be
W({\bf s}^i - x_1 {\bf c}^i_1,\bm D'_{1i}),
\quad
\bm D'_{1i} = \{{\bf c}^i_2,{\bf c}^i_3,\ldots,{\bf c}^i_m\},
\label{Zero13a}
\ee
{\it i.e.}, the matrix $\bm D'_{1i}$ is obtained from the matrix ${\bf D}_i$ by dropping the 
first column ${\bf c}^i_1$.

The definition (\ref{WvectGF}) of the generating function for the vector partition $W({\bf s},{\bf D})$ 
implies the recursion
\be
W({\bf s},{\bf D}) - W({\bf s}-{\bf c}_1,{\bf D}) = W({\bf s},\{c_{2},c_{3},\ldots,c_{m}\}),
\label{VPFrecursion}
\ee
leading to
\be
W({\bf s}^i ,{\bf D}_i) - W({\bf s}^i - {\bf c}^i_1,{\bf D}_i) = W({\bf s}^i,\bm D'_{1i}).
\label{VPFrecursion1}
\ee
Use this relation to estimate the infinite sum of the terms (\ref{Zero13a})
to obtain
\be
\sum_{x_1=0}^{\infty} W({\bf s}^i - x_1 {\bf c}^i_1,\bm D'_{1i})  =
\sum_{x_1=0}^{\infty} [W({\bf s}^i - x_1 {\bf c}^i_1,{\bf D}_i) - W({\bf s}^i - (x_1+1) {\bf c}^i_1,{\bf D}_i)] = 
W({\bf s}^i ,{\bf D}_i),
\label{Zero14}
\ee
as there always exists a positive integer $k$ such that $W({\bf s}^i - k {\bf c}^i_1,{\bf D}_i)$
evaluates to zero. It is true if at least one element of ${\bf c}^i_1 = c_{il}\hat{\bf c}_1$
is positive 
which is equivalent to the existence of at least one positive element in the column ${\bf c}_1$.
If this condition fails, we apply the column sign inversion procedure (\ref{mCayley118}) to the 
column ${\bf c}_1$ with the corresponding update of the column ${\bf s} \Rightarrow {\bf s} - {\bf c}_1$.
Write the vector partition $W({\bf s},{\bf D})$ as 
$$
W({\bf s},{\bf D})  = C_{\bf s}^1 + \sum_{i=2}^m W({\bf s}^i,{\bf D}_i),
$$
where $C_{\bf s}^1$ denotes (yet unknown) contribution of the column ${\bf c}_1$. 
The relation (\ref{Zero14}) shows that each contribution $W({\bf s}^i,{\bf D}_i)$ for $2 \le i \le m$
can be written as the infinite sum of $W({\bf s}^i - x_1 {\bf c}^i_1,\bm D'_{1i})$. 
This leads to 
\be 
W({\bf s},{\bf D}) =\sum_{x_1=0}^{\infty} W({\bf s}'_1,\bm D'_1) = 
\sum_{i=2}^m \sum_{x_1=0}^{\infty} W({\bf s}^i - x_1 {\bf c}^i_1,\bm D'_{1i}) =
\sum_{i=2}^m W({\bf s}^i,{\bf D}_i),
\label{Zero15}
\ee
and the contribution $C_{\bf s}^1$ vanishes.
This result is used in the example discussed in Appendix \ref{Triple2}.


If the matrix $\bm D'_{1}$ has another LZ-column ${\bf c}_j$ 
we have two options -- either the columns ${\bf c}_1$ and ${\bf c}_j$
are collinear or they are not.  
In the case of noncollinear LZ-columns the contributions $C_{\bf s}^i$ are independent 
of each other and thus we have $C_{\bf s}^1=C_{\bf s}^j=0$. 
Application of the Cayley reduction algorithm (or its extended version) to the 
partition $W({\bf s},{\bf D})$ with two collinear LZ-columns ${\bf c}_i$ and ${\bf c}_j$ produces
two partitions $W({\bf s}^i ,{\bf D}_i)$ and $W({\bf s}^j ,{\bf D}_j)$ where 
both matrices ${\bf D}_i, {\bf D}_j$ have one column made of only zeros. 
Generally such a partition would correspond to the infinite number of solutions of the 
system ${\bf s}^i = {\bf D}_i \cdot {\bf x}$, but it is not true in this specific case. 
As the vector ${\bf s}^i$ has at least one {\it negative} element
the corresponding system has no solutions and thus $C_{\bf s}^i=C_{\bf s}^j=0$.



\section{Vector partitions with C-columns}
\label{Collinear}

The reduction algorithms 
described in Sections \ref{Cayley},\ref{Multi} cannot be applied in case when
in the generator matrix ${\bf D}$ has the collinear columns (C-columns) as 
it leads to appearance of the zero columns in the reduced partition.
There are two ways to overcome this problem. The first one considered 
in Section \ref{Convolution} is quite 
straightforward and it just rewrites the system ${\bf s}={\bf D} \cdot {\bf x}$ 
transferring all C-columns to the l.h.s. so that the modified 
matrix 
can be reduced using the algorithms discussed in the 
preceding Sections. 

The second approach is more elaborated as it is based on the extension 
of the original equations ${\bf s}={\bf D} \cdot {\bf x}$ to a
larger but equivalent system $\tilde {\bf s}=\tilde{\bf D} \cdot \tilde {\bf x}$
where the matrix $\tilde {\bf D}$ is void of C- and NP-columns.
The extended partition can be then reduced by the Cayley elimination procedure 
and we consider the details of this process in Section \ref{Extension}.
Finally, in Section \ref{Comparison} we compare these two
approaches.


\subsection{Partition convolution superposition}
\label{Convolution}

Call an integer generator matrix {\it competent} if its columns (viewed as vectors) are noncollinear
and {\it incompetent} otherwise.
Assume that the vectors corresponding to the first $n$ columns of the incompetent generator 
matrix ${\bf D}$ are parallel to each other
and rewrite the linear system ${\bf s}={\bf D} \cdot {\bf x}$ as follows
\be
 {\bf s} = \sum_{i=1}^m {\bf c}_i x_i,
\quad
{\bf c}_j = u_j {\bf C}, 
\quad
{\bf C} = \{C_1,C_2,\ldots,C_l\}^T,
\quad
1 \le j \le n < m,
\label{noncol01}
\ee
where $u_j$ are positive integers and one can always select the vector ${\bf C}$ 
to be a P-column.

The problem (\ref{noncol01}) is  equivalent to 
\be
{\bf S}(k) = {\bf s}-k {\bf C} = \sum_{i=n+1}^m {\bf c}_i x_i,
\quad
k = \sum_{i=1}^n u_i x_i, 
\quad
0 \le k \le K = \min(\lfloor s_r/C_r \rfloor ),
\quad
1 \le r \le l,
\label{noncol02}
\ee
where $K$ is computed using only {\it positive} $C_r > 0$.
Introduce a vector ${\bf u} = \{u_1,u_2,\ldots,u_n\}$ and a competent matrix with $(m-n)$ columns
$\bar{\bf D} = \{{\bf c}_{n+1},\ldots,{\bf c}_m\},$
for which the corresponding partition $W({\bf S}(k),\bar{\bf D})$ 
admits the reduction either as the sum of $W_i^l$ or the mixture of $W^l_i$ and $\bar W^l_i$.
The number of solutions of the system (\ref{noncol02}) for given value of $k$  is equal to 
a number of solutions of the first equation given by the partition $W({\bf S}(k),\bar{\bf D})$
multiplied by a number of solutions of the second equation, {\it i.e.}, by a scalar partition $W(k,{\bf u})$.
Then the vector partition $W({\bf s},{\bf D})$ is equivalent 
to a convolution
\be
W({\bf s},{\bf D}) = 
\sum_{k=0}^{K} W(k,{\bf u})H({\bf s}-k {\bf C})W({\bf s}-k {\bf C},\bar{\bf D}),
\quad
K = \min(\lfloor s_r/C_r \rfloor ),
\quad
1 \le r \le l,
\label{noncol04}
\ee 
of the scalar partition and the vector partition $W({\bf s}-k {\bf C},\bar{\bf D})$ built on the competent matrix $\bar{\bf D}$
(it admits the reduction described in Section \ref{Multi}).
Note that the function $H({\bf s}-k {\bf C})$ in (\ref{noncol04}) is applied 
only to the elements of the vector ${\bf s}-k {\bf C}$ corresponding to the rows
of $\bar{\bf D}$ void of negative generators.

This approach can be extended to the case when the matrix ${\bf D}$ has $r > 1$ groups of the 
C-columns. Each $p$-th group is characterized by the P-column ${\bf C}_p$
and the $n_p$ integers 
$u_{p,i}$ making up the vector  ${\bf u}_p = \{u_{p,1},u_{p,2},\ldots,u_{p,n_p}\}$. 
Introduce the vector ${\bf k} =\{k_{1},k_{2},\ldots,k_{r}\}$ and  (\ref{noncol02}) turns into 
\be
{\bf S}({\bf k}) = {\bf s}-\sum_{p=1}^r k_p {\bf C}_p = \sum_{i=n+1}^m {\bf c}_i x_i,
\quad
k_p = \sum_{i=1}^{n_p} u_{p,i} x_i, 
\quad
n = \sum_{p=1}^r n_p.
\label{noncol05}
\ee
The vector partition $W({\bf s},{\bf D})$ evaluates to the 
convolution
\be
W({\bf s},{\bf D})= 
\sum_{\bf k=0}^{\bf K} 
H({\bf S}({\bf k}))W({\bf S}({\bf k}),\bar{\bf D})  \cdot  \prod_{p=1}^r W(k_p,{\bf u}_p),
\
K_p =  \min(\lfloor s_i/C_{p,i}\rfloor ),
\
\sum_{\bf k=0}^{\bf K}\equiv \sum_{k_1=0}^{K_1} 
\ldots \sum_{k_r=0}^{K_r}.
\label{noncol06}
\ee 
Consider a special case $n=m-1$, {\it i.e.}, when the competent matrix $\bar{\bf D}=\{{\bf c}\}$ has a single column ${\bf c}$.
The vector partition $W({\bf S}({\bf k}),\{{\bf c}\})=1$ 
only when ${\bf S}({\bf k}) = k_0 {\bf c}$ for some nonnegative integer $k_0$,
otherwise it vanishes. Thus we observe that
$W({\bf S}({\bf k}),\bar{\bf D}) = \delta_{{\bf S}({\bf k}),k_0 {\bf c}}$ 
and we obtain using (\ref{noncol06}) 
\be
W({\bf s},{\bf D})= 
\sum_{\bf k=0}^{\bf K} \sum_{k_0=0}^{K_0} 
\delta_{{\bf S}({\bf k}), k_0 {\bf c}}H({\bf S}({\bf k})) \cdot  \prod_{p=1}^r W(k_p,{\bf u}_p),
\quad
K_0 = \max(\lfloor s_i/c_i \rfloor ), 
\label{noncol07}
\ee 
where again only positive $c_i$ are used for $K_0$ evaluation.

The extreme case 
arises for $n=m$, 
when {\it each} column of the matrix ${\bf D}$ belongs to one of $r$  groups of the C-columns.
Use ${\bf k}' =\{k_{1},k_{2},\ldots,k_{r-1}\}$ to write
\be
{\bf S}({\bf k}') = {\bf s}-\sum_{p=1}^{r-1} k_p {\bf C}_p = k_r {\bf C}_r = \sum_{i=1}^{n_r} u_{r,i} x_i {\bf C}_r .
\label{noncol08}
\ee
As ${\bf S}({\bf k}')$ and ${\bf C}_r$ are  collinear we obtain a set    
of additional conditions on $k_r$ leading to
\be
W({\bf s},{\bf D})= 
\sum_{{\bf k} =0}^{{\bf K}} 
\delta_{{\bf S}({\bf k}'),k_r {\bf C}_r}H({\bf S}({\bf k}')) \cdot  \prod_{p=1}^r W(k_p,{\bf u}_p),
\quad
K_r =  \max(\lfloor s_i/C_{r,i} \rfloor ).
\label{noncol09}
\ee 
The example of partitions with C-columns is discussed in Appendix \ref{CollinearExample1}
where we describe the reduction procedure for the matrix used in computation of the 
Kronecker coefficients \cite{Mishna2021,Mishna2024}.

\subsection{Collinear columns elimination}
\label{Extension}

An alternative approach to the reduction of the incompetent matrix 
was introduced by Sylvester in the fifth lecture of his course on the partitions \cite{SylvLectures1859}
where he suggested to supply the system ${\bf s} = {\bf D} \cdot {\bf x}$
with a new unknown $x_{0}$ and a single auxiliary equation to generate a new system   
$\tilde {\bf s} = \tilde {\bf D} \cdot \tilde {\bf x}$
with $\tilde {\bf x} = \{x_0,x_1,x_2,\ldots,x_m\}$.
The $(l+1)\times (m+1)$ integer competent matrix 
$\tilde {\bf D}=\{\tilde {\bf c}_{0},\tilde {\bf c}_1,\tilde {\bf c}_2,\ldots,\tilde {\bf c}_m\}$ 
consists of $(m+1)$ columns $\tilde {\bf c}_{i} =\{c_{i0},c_{i1},c_{i2},\ldots,c_{il}\}^T$.
The elements $c_{i0}$ and the value of $s_0$ in the added equation must be chosen to satisfy two conditions -- 
(a) each nonnegative integer solution of the original system ${\bf s} ={\bf D} \cdot {\bf x}$
should correspond to a single solution of the system $\tilde {\bf s} = \tilde {\bf D} \cdot \tilde {\bf x}$,
and (b) the matrix  $\tilde {\bf D}$ should be void of C- and NP-columns.
The first condition is equivalent to 
\be
W(\tilde{\bf s},\tilde{\bf D}) = W({\bf s},{\bf D}),
\label{solution_number}
\ee
while the second one guarantees that the vector partition $W(\tilde{\bf s},\tilde{\bf D})$ 
admits the first step of variable elimination procedure of reduction described in Section \ref{Cayley}.
This approach was recently employed for derivation 
of a new class of linear relations involving scalar partitions \cite{RubLinRel2025}.

The choice of the added equation is not unique and following the 
Sylvester advice we use
$$
s_0 = x_{0} - \sum_{i=1}^{m} \delta_i x_i, 
$$
with positive $\delta_i > 0$ and nonnegative $s_0$.
Thus $c_{i0} = -\delta_i$ and the resulting system
corresponds to the vector partition $W(\tilde{\bf s},\tilde{\bf D})$ with 
\be
\tilde{\bf D} =
\left [
\begin{array}{ccccc}
1 &-\delta_1 & -\delta_2 & \ldots & -\delta_m   \\
0 & c_{11} & c_{21} & \ldots & c_{m1} \\
0 &c_{12} & c_{22} & \ldots & c_{m2} \\
\ldots   & \ldots  & \ldots  & \ldots \\
0 &c_{1l} & c_{2l} & \ldots & c_{ml} \\
\end{array}
\right ],
\quad
\quad
\tilde{\bf s} = \left [
\begin{array}{c}
s_0 \\
s_1 \\
s_2 \\
\ldots \\
s_m
\end{array}
\right ].
\label{matr_hatD}
\ee
Each integer nonnegative solution 
${\bf x}=\bm\xi =\{\xi_1,\xi_2,\ldots,\xi_m\}$
of the system ${\bf s} ={\bf D} \cdot {\bf x}$
corresponds to the following solution $\tilde {\bf x}$ of the extended system
\be
x_i = \xi_i, 
\
(1 \le i \le m),
\quad\quad
x_{0} = 
\sum_{i=1}^{m} \delta_i \xi_i,
\label{system2sol}
\ee
where we set $s_0=0$. 
We observe that the condition (\ref{solution_number}) is satisfied as for each solution
$\bm\xi$ of the original Diophantine equation $s = {\bf d} \cdot {\bf x}$ 
the system $\tilde {\bf s} = \tilde {\bf D} \cdot \tilde {\bf x}$ has
exactly one solution (\ref{system2sol}).

The selection of the positive $\delta_i$ values is a (partially heuristic) process
governed by the requirement that all columns of the matrix $\tilde {\bf D}$ should 
have $\mbox{gcd}(\tilde{\bf c}_i)=1$ and be 
not collinear to each other.
To satisfy the first condition one has to choose $\delta_i$  
that does not divide the GCD of the column ${\bf c}_i$ of the original matrix ${\bf D}$, i.e., $\delta_i \nmid \mbox{gcd}({\bf c}_i)$.
As for the second condition one has to consider several cases. 
Start from the matrix ${\bf D}$ having a single group of  $n$
C-columns ${\bf c}_j = u_j {\bf C}$ with integer $u_j$
and ${\bf C}$ chosen to be a P-column.
To make the 
$\tilde {\bf c}_j$ noncollinear P-columns use the 
following simple rules 
\be
\mbox{gcd}(u_j,\delta_j)=1,
\quad
\delta_j/u_j \ne \delta_k/u_k,
\quad 
1 \le j \ne k \le n.
\label{one_col_group}
\ee
When the matrix ${\bf D}$ has several groups of the 
C-columns, the columns belonging to the different groups are 
not collinear to each other and thus the procedure of $\delta_i$ selection
can be done using
the rules (\ref{one_col_group}) independently within each group. 
The search of $\delta_i$ is discussed in 
the illustrative example presented in Appendix \ref{CollinearExample2}.

It should be underlined that the rules (\ref{one_col_group}) leave a lot of freedom in the 
$\delta_i$ selection so that the matrix $\tilde {\bf D}$ used for the further reduction 
is not unique.
The reduction itself is done by using the Cayley elimination algorithm that 
produces $(m+1)$ vector partitions $W(\tilde{\bf s}^i,\tilde{\bf D}^i)$, but 
$W(\tilde{\bf s}^{0},\tilde{\bf D}^{0})$ evaluates to zero as explained in 
Section \ref{ZeroMulti} and we end up with just $m$ vector partitions. 

From comparison of the matrices ${\bf D}$ and $\tilde {\bf D}$
it follows that the elements of the vector ${\bf s}^1$ given by (\ref{Gen407}) and the components 
of $\tilde{\bf s}^1$ from the second to the last one coincide -- $s^1_i = \tilde{s}^1_{i+1},\ 1\le i \le l$.
The same is true for the columns ${\bf c}^1_i$  and $\tilde{\bf c}^1_i$ 
with $2 \le i \le m$. 
We already noted in Section \ref{Multi} for 
the matrix ${\bf D}$ with two C-columns ${\bf c}_1 = u_1 {\bf C}$
and ${\bf c}_k = u_k {\bf C}$ after reduction we obtain for
the column ${\bf c}^1_k = {\bf 0}$. This means that 
in the column $\tilde{\bf c}^1_k$ the first element 
$\tilde{c}^1_{k1} = (u_k \delta_1 - u_1 \delta_k) C_m$ is nonzero while all other elements vanish,
and the contribution of $\tilde{\bf c}^1_k$ in the next reduction round equals zero.

\subsection{Comparison of the methods}
\label{Comparison}

In the convolution method in Section \ref{Convolution}
the transfer of a group of the C-columns from the matrix ${\bf D}$ to the argument ${\bf S}$
diminishes the number of columns in $\bar {\bf D}$ so that it may significantly 
decrease the number of the produced scalar partitions. On the other hand, 
this relatively small number of the scalar partitions enter the 
summands in the multiple sums that for the large values of ${\bf s}$ may lead to 
very large number of terms requiring evaluation. 

Turning to the alternative matrix extension method in Section \ref{Extension}
we note that it replaces a single $l \times m$ matrix by $m$ matrices of the same size
that might significantly increase 
the total number of the scalar partitions in the final expression. 
The same time this larger number of scalar partitions are evaluated only once 
to obtain the value of the original vector partition. 
For a given vector partition this algorithm also allows to generate multiple {\it equivalent}
superpositions of scalar partitions.


\section{Iterative reduction of vector to scalar partitions}
\label{Muirhead}
The basic step of iterative reduction procedure
of the vector partition corresponding to the 
$l$-row matrix ${\bf E}$ described in Section \ref{Multi} 
includes elimination of the last row and a single $i$-th column of 
${\bf E}$. It results in a single partition based on the $(l-1)$-row
matrix 
for the P-column or a superposition of 
partitions built on similar matrices 
for the NP-column. Each such matrix undergoes 
similar transformation in the consecutive iterative steps leading to a combination of scalar partitions.
The total number of iterations required to reach this stage is $l-1$.

We note in Section \ref{Cayleyorig} that the order of the rows in the 
two-row matrix ${\bf E}$ affects the individual scalar partition contributions 
to the double partition. The same holds for a matrix with the larger number of rows.
This means that while one can change row order in the initial matrix ${\bf E}$,
it is forbidden for any matrix generated in the course of the iterative elimination procedure. 
For definiteness assume that the order of row elimination 
is fixed and at every step we eliminate the last row of each generated matrix.

Consider a vector partition $A_{\nu} W({\bf E}_{\nu})$ corresponding to the
$n \times (p+1)$ matrix ${\bf E}_{\nu}$ obtained after $\nu$ reduction steps
of the original 
$l \times (m+1)$ matrix ${\bf E}$, so that $n=l-\nu$ and $p=n-\nu$. The factor $A_{\nu}$ 
is a product of the factors shown in (\ref{Gen407lastrowpositive}) and (\ref{Gen411GCDlastrowpositive}) 
computed at each reduction step.
The matrix ${\bf E}_{\nu}$ can be presented as a set of rows ${\bf E}_{\nu} = \{{\bf R}_1,{\bf R}_2,\ldots,{\bf R}_n\}$, where 
the row ${\bf R}_k=\{S_k=R_{k0},R_{k1},\ldots,R_{kp}\}$ has $p+1$  elements and 
$S_k$ is the linear superposition of the components $s_i$ of the original 
$l \times (m+1)$ matrix ${\bf E}$
$$
S_k = L'_k({\bf s}) = \kappa_{k0} +  \sum_{j=1}^m \kappa_{kj} s_j,
$$
with integer coefficients $\kappa_{kj}$.
The last row ${\bf R}_n=\{S_n=R_{n0},R_{n1},\ldots,R_{np}\}$ has all nonnegative numerical elements   
while other rows ${\bf R}_k\ (1 \le k \le n-1)$ might have negative elements.

How can one define the partition function in this case? 
The last row ${\bf R}_n$ of the matrix corresponds to a scalar partition
defined by (\ref{coin1}) that requires to find all sets ${\bf x}=\bm \xi_K$ of nonnegative integer solutions.
Assume that we find all such sets, 
their number is given by a scalar partition 
$W(S_n,\{R_{n1},\ldots,R_{np}\})$. Now substitute each solution $\bm \xi_K$ into the 
remaining $n-1$ equations and retain only those that satisfy them all. The number of such solutions
is the vector partition $W({\bf E}_{\nu})$.
This definition fails when the last row has zero entries 
as the number of solutions $\bm \xi_K$ in this case is infinite. 
Nevertheless, one still can proceed with the evaluation procedure in these cases,
this problem is addressed in Section \ref{ZeroMulti} where we prove that the 
contribution of a column with the last zero element is zero. 

In the original definition of the vector partition mentioned in Section \ref{intro0} the 
elements of the vector ${\bf s}$ are considered nonnegative. This assumption does not hold in general 
for ${\bf E}_{\nu}$ as the  linear superpositions $L'_k({\bf s})$ might be negative for some values of the 
vector ${\bf s}$. For the rows with negative numerical elements the negative $L'_k({\bf s}) < 0$ are acceptable.
But if the row ${\bf R}_k$ is void of negative numbers
we have to restrict its initial element $L'_k({\bf s})$ to the nonnegative values $L'_k({\bf s}) \ge 0$ 
and to replace the partition $A_{\nu} W({\bf E}_{\nu})$ by
\be
A_{\nu} W({\bf E}_{\nu}) \Rightarrow A_{\nu} B_{\nu} W({\bf E}_{\nu}), 
\quad
B_{\nu} = \prod_{k} H(L'_k({\bf s}) ).
\label{UnitStep2c}
\ee
Application of the procedure (\ref{UnitStep2c}) after every 
reduction step is required for determination of the partition chambers.


\section{Conclusive remarks}
\label{Discuss}

In this manuscript we present an implementation of the Sylvester program aimed to reduction of
vector (compound) partition to scalar partitions  \cite{Sylv1,SylvLectures1859}. The first step in this direction
was made in 1860 by Cayley \cite{Cayley1860} when he showed that
the double partition problem subject to a set of conditions
admits a compact solution
in which the $i$-th column ${\bf c}_i$ of a {\it positive} generator matrix  ${\bf D}$ 
leads to a single scalar partition contribution $W_i^2$ (see Section \ref{Cayleyorig}).
The author of this manuscript have suggested
an alternative approach that solves the problem
when the original Cayley method fails because the conditions are not satisfied \cite{RubDouble2023} .

The results presented in the manuscript indicate that Sylvester strongly 
underestimated the number of scalar partitions required for 
computation of the vector partition --
their minimal number in case of $l \times m$ positive matrix might significantly surpass
$m!/(m-l+1)!$ which in its turn is $(l-1)!$ times larger than predicted by Sylvester in \cite{Sylv1}.
We address the case of collinear columns
in the matrix ${\bf D}$ in Section \ref{Collinear} and show that 
a vector partition either boils down to a convolution of 
scalar partitions or it is replaced by an equivalent one with the
larger matrix having no such columns.

Our approach is a generalization of the algorithms discussed in \cite{RubDouble2023}
and here we demonstrate that a vector partition always 
reduces iteratively  either to a sum of scalar partitions themselves or their convolutions,
and thus  the Sylvester program 
is successfully implemented.
Recalling that the scalar partitions can be computed as 
a superposition of the known polynomials \cite{Rub04}
we arrive at the important statement --- {\it any vector partition} 
can be expressed through the Bernoulli polynomials of higher order.

All known approaches to evaluation of vector partitions 
(see \cite{Barvinok1993, Sturmfels1995, Brion1997, Szenes2003}), their later modifications 
\cite{Beck2004, DeLoera2004, Milev2023}
and the software packages based of these methods successfully determine
the chambers but fail to produce the explicit expressions for the vector partitions as they 
replace the original generating function by the signed sum of the simpler partial generating functions
valid in specific chambers.
This means that for each new argument column ${\bf s}$ one has to 
perform the Taylor expansion of the partial generating functions to extract the 
desired coefficient. 
Our method is the first one that permits the derivation of a {\it single explicit formula} for the 
vector partition which is valid for an {\it arbitrary} value of ${\bf s}$. The usage of this formula does not 
require the explicit knowledge of the chambers but nevertheless it implicitly contains their description.




\newpage
{\LARGE \bf Appendices}

\appendix
\section{Triple partition reduction example}
\renewcommand{\theequation}{\thesection\arabic{equation}}
\setcounter{equation}{0}
\label{Triple1}

Consider a reduction of a triple partition $W({\bf E})$ defined by a matrix ${\bf E}$
\be
{\bf E} = \left[
\begin{array}{ccccc}
s_1 & 1 & 2 & 1 & 1 \\
s_2 & 2 & 1 & 3 & 4 \\
s_3 & 3 & 2 & 3 & 2
\end{array}
\right],
\label{A1}
\ee
and perform the first reduction step by eliminating the last row of ${\bf E}$.
We obtain 
\be
W({\bf E}) = \sum_{k=1}^4 W({\bf E}_k), 
\label{A2}
\ee
with 
\bea
&
{\bf E}_1 = \left[
\begin{array}{cccc}
3s_1-s_3 & 4 & 0 & 1 \\
3s_2-2s_3 & -1 & 3 & 8
\end{array}
\right],
\quad
{\bf E}_2 = \left[
\begin{array}{cccc}
3s_1-s_3 & 0 & 4 & 1 \\
3s_2-3s_3 & -3 & -3 & 6
\end{array}
\right],
\nonumber
\\
&
{\bf E}_3 = \left[
\begin{array}{cccc}
2s_1-s_3 & -1 & 2 & -1 \\
2s_2-4s_3 & -8 & -6 & -6
\end{array}
\right],
\quad
{\bf E}_4 = \left[
\begin{array}{cccc}
2s_1-2s_3 & -4 & -4 & -2 \\
2s_2-s_3 & 1 & 3 & 6
\end{array}
\right],
\label{A2a}
\eea
where each matrix ${\bf E}_i$ contains the NP-column.
Perform the following transformation of the matrices ${\bf E}_i$ -- first divide each row by 
the GCD of its numerical elements, and then use (\ref{mCayley118}, \ref{Gen407lastrowpositive})
to make last row positive, and obtain a sum of four double partitions
\be
W({\bf E}) = \sum_{k=1}^4 a_k W(\tilde{\bf E}_k), 
\quad
a_1, a_3 = -1, 
\quad a_2, a_4 = 1,
\label{A3}
\ee
with
\bea
&
\tilde{\bf E}_1 = \left[
\begin{array}{cccc}
3s_1-s_3+4 & -4 & 0 & 1 \\
3s_2-2s_3-1 & 1 & 3 & 8
\end{array}
\right],
&
\tilde{\bf E}_2 = \left[
\begin{array}{cccc}
3s_1-s_3+4 & 0 & -4 & 1 \\
s_2-s_3-2 & 1 & 1 & 2
\end{array}
\right], 
\nonumber
\\
&
\tilde{\bf E}_3 = \left[
\begin{array}{cccc}
2s_1-s_3 & 1 & -2 & 1 \\
s_2-2s_3-10 & 4 & 3 & 3
\end{array}
\right],
&
\tilde{\bf E}_4 = \left[
\begin{array}{cccc}
s_1-s_3 & -2 & -2 & -1 \\
2s_2-s_3 & 1 & 3 & 6
\end{array}
\right].
\label{A4}
\eea
As discussed in Section \ref{Muirhead}   
each summand in (\ref{A3}) acquires an additional factor
producing 
\be
W({\bf E}) = \sum_{k=1}^4 A_k W(\tilde{\bf E}_k), 
\label{A5}
\ee 
where 
\be
A_1 = -
H(3s_2-2s_3-1), \quad
A_2 =
H(s_2-s_3-2), \quad
A_3 = -
H(s_2-2s_3-10), \quad
A_4 = 
H(2s_2-s_3).
\label{A6}
\ee
Each double partition in (\ref{A5}) to be reduced to a set of scalar partitions
as shown in Section \ref{Cayley}. 
Note that $\tilde{\bf E}_1$ 
has a single NP-column $\{0,3\}^T$ and
thus, the elimination of this column should be done using the general Cayley algorithm 
presented in Section \ref{Cayleynew}, while all other columns are processed as shown in 
Section \ref{Cayleyorig}.
The results for $W(\tilde{\bf E}_k)$ 
read
\bea
W(\tilde{\bf E}_1) &=& 
W(s_1+4s_2-3s_3,\{4,11\}) + W(8s_1-s_2-2s_3-1,\{1,11\}) 
\nonumber \\
&-& \sum_{j=2}^4   W((3s_1-s_3+4)/3-j,\{1,4\}) \delta_{(3s_2-2s_3-1) \bmod{3},0}
\nonumber \\
& - &  \sum_{j=1,2,4}  W((3s_1-s_3+2)/3-j,\{1,4\}) \delta_{(3s_2-2s_3-2) \bmod{3},0}
\nonumber \\
& - &\sum_{j=0,2,3} W((3s_1-s_3)/3-j,\{1,4\}) \delta_{(3s_2-2s_3) \bmod{3},0}\;,
\nonumber 
\eea
\bea
W(\tilde{\bf E}_2) &=&
W(3 s_1 + 4 s_2 - 5 s_3-4, \{9, 4\})
- W(3 s_1 - s_3, \{1, 4\})
+ W(6 s_1 - s_2 - s_3-4, \{1,9\}),
\nonumber \\
W(\tilde{\bf E}_3)  &=&
W(6 s_1 + 2 s_2 - 7 s_3-20, \{9, 11\})
- W(8 s_1 -s_2 - 2s_3-1, \{1, 11\})
\nonumber \\
&+& W(6 s_1 - s_2 - s_3, \{1,9\}),
\nonumber \\
W(\tilde{\bf E}_4)  &=&
W(s_1 + 4 s_2 - 3 s_3, \{4,11\})
- W(3 s_1 + 4 s_2 - 5 s_3-4, \{4,9\})
\nonumber \\
&+& W(6 s_1 + 2 s_2 - 7 s_3-20, \{9,11\}). 
\nonumber
\eea
These expressions together with (\ref{A5},\ref{A6}) give the 
complete reduction of the triple partition defined by the matrix (\ref{A1}) into
a sum of twenty scalar partitions. 
Note that the first argument $L({\bf s})$ of each scalar partition 
should be nonnegative that together with the factor $A_i$ in (\ref{A5}) 
completely determine the  scalar partition chamber.

\section{Alternative triple partition reduction}
\renewcommand{\theequation}{\thesection\arabic{equation}}
\setcounter{equation}{0}
\label{Triple2}

In the Appendix \ref{Triple1} after the first reduction step we encounter
a matrix $\tilde{\bf E}_1$ having NP-column $\{0,3\}^T$ that requires 
application of the general Cayley algorithm for its reduction. 
The same time in Section \ref{ZeroMulti} we prove that the reduction of a column of the type $\{b,0\}^T$ with positive $b$
leads to vanishing contribution. It is instructive to reduce the number of terms in the 
final expression of the triple partition via the scalar ones. 
The convenient way to do that is 
to reshuffle the rows in the original triple matrix ${\bf E}$.

Consider reduction of a triple partition $W({\bf E}')$ defined by the matrix
\be
{\bf E}' = \left[
\begin{array}{ccccc}
s_2 & 2 & 1 & 3 & 4 \\
s_1 & 1 & 2 & 1 & 1 \\
s_3 & 3 & 2 & 3 & 2
\end{array}
\right],
\label{B1}
\ee
obtained from ${\bf E}$ by reordering of the first two rows.
We find
\be
W({\bf E}) = W({\bf E}') = \sum_{k=1}^4 a'_k W(\tilde{\bf E}'_k), 
\quad
a'_1, a'_3 = -1, 
\quad a'_2, a'_4 = 1,
\label{B3}
\ee
where after the required transformations we obtain
\bea
&
\tilde{\bf E}'_1 = \left[
\begin{array}{cccc}
2s_2-s_3 & -1 & -3 & -6 \\
s_1-s_3-5 & 2 & 2 & 1 
\end{array}
\right],
\quad
\tilde{\bf E}'_2 = \left[
\begin{array}{cccc}
3s_2-2s_3 & -1 & 8 & 3 \\
3s_1-s_3 & 4 & 1 & 0 
\end{array}
\right],
\nonumber
\\
&
\tilde{\bf E}'_3 = \left[
\begin{array}{cccc}
s_2-s_3-1 & 1 & -1 & 2 \\
3s_1-s_3 & 0 & 4 & 1 
\end{array}
\right],
\quad
\tilde{\bf E}'_4 = \left[
\begin{array}{cccc}
s_2-2s_3-7 & -3 & 4 & 3 \\
2s_1-s_3-2 & 2 & 1 & 1
\end{array}
\right].
\label{B4}
\eea
Adding the $H$ factors we have
\begin{equation*}
W({\bf E}') = \sum_{k=1}^4 A'_k W(\tilde{\bf E}'_k), 
\
A'_1 = -
H(s_1-s_3-5), \
A'_2 = -A'_3 =
H(3s_1-s_3), \
A'_4 = 
H(2s_1-s_3-2).
\end{equation*}
The second reduction step turns each double partition $W(\tilde{\bf E}'_k)$
into a sum of scalar partitions, but this time the contribution of the LZ-columns $\{3,0\}^T$ in $\tilde{\bf E}'_2$ 
and $\{1,0\}^T$ in $\tilde{\bf E}'_3$ vanishes
and we find 
\bea
W(\tilde{\bf E}'_1) &=&
W(s_1 + 4 s_2 - 3 s_3, \{4,11\})
- W(3 s_1 + 4 s_2 - 5 s_3-4, \{4,9\})\nonumber \\
&+&
W(6 s_1 + 2 s_2 - 7 s_3-20, \{9,11\}),
\label{B7} \\
W(\tilde{\bf E}'_2)  &=&
W(3 s_1 + 4 s_2 - 5 s_3, \{11, 4\})
- W(-8 s_1 + s_2 + 2 s_3-11, \{1,11\}),
\nonumber \\
W(\tilde{\bf E}'_3)  &=&
W(3 s_1 + 4 s_2 - 5 s_3, \{9, 4\})
- W(-6 s_1 + s_2 + s_3-10, \{1,9\}),
\nonumber \\
W(\tilde{\bf E}'_4)  &=&
W(6 s_1 + 2 s_2 - 7 s_3-20, \{9, 11\})
+ W(-8 s_1 +s_2 + 2s_3-11, \{1, 11\})
\nonumber \\
&-&
W(-6 s_1 + s_2 + s_3-10, \{1,9\}).
\nonumber 
\eea
The total number of scalar partitions in (\ref{B7}) is ten which is 
significantly smaller than we have
in Appendix \ref{Triple1}.

\section{Partition with collinear columns -- convolution algorithm}
\renewcommand{\theequation}{\thesection\arabic{equation}}
\setcounter{equation}{0}
\label{CollinearExample1}

Consider reduction of a triple partition $W({\bf E})$ defined by the matrix
used in the computation of the Kronecker coefficients \cite{Mishna2021,Mishna2024}
\be
{\bf E} = \left[
\begin{array}{cccccccccccc}
s_3 & 0 & 0 & 1 & 1 & 1 & 1 & 2 & 1 & 2 & 2 & 3  \\
s_2 & 0 & 1 & 0 & 0 & 1 & 1 & 1 & 1 & 1 & 2 & 2 \\
s_1 & 1 & 0 & 0 & 1 & 0 & 0 & 0 & 1 & 1 & 1 & 1 
\end{array}
\right].
\label{C1}
\ee
It has two equal LZ-columns $\{1,1,0\}^T$ but their 
contribution vanishes 
as discussed in Section \ref{ZeroMulti}. 
As the last row of ${\bf E}$ has five zeros the  
first reduction step produces only six nonvanishing terms
\bea
&W({\bf E}) = \sum_{k=1}^6 A_6 W({\bf E}_k), \quad
A_1 = H(s_2)H(s_3), &
A_2 = -H(s_2-2s_1-6)H(s_3-3s_1-9), 
\label{C2}\\
&
A_3 = -H(s_2)H(s_3-s_1-1), &
A_4 = H(s_2-2s_1-6)H(s_3-2s_1-4) ,
\nonumber  \\
&
A_5 = H(s_2-s_1-2)H(s_3-s_1-1) ,&
A_6 = -H(s_2-s_1-2)H(s_3-2s_1-4),
\nonumber  
\eea
and the matrices ${\bf E}_k$ read
\bea
&
{\bf E}_1 = \left[
\begin{array}{ccccc|cc|cc|cc}
s_3 & 1 & 1 & 1 & 2 & 2 & 2 & 1 & 1 & 3 & 0  \\
s_2 & 1 & 1 & 1 & 2 & 1 & 1 & 0 & 0 & 2 & 1 
\end{array}
\right],
\nonumber
\\
&
{\bf E}_2 = \left[
\begin{array}{ccccc|cc|cc|cc}
s_3-3s_1-9 & 1 & 1 & 1 & 2 & 2 & 2 & 1 & 1 & 3 & 0  \\
s_2-2s_1-6 & 1 & 1 & 1 & 2 & 1 & 1 & 0 & 0 & 2 & 1 
\end{array}
\right],
\nonumber
\\
&
{\bf E}_3 = \left[
\begin{array}{ccccc|cc|cc|cc}
s_3-s_1-1 & 1 & 1 & 1 & 2 & 0 & 0 & 1 & 1 & 1 & 2  \\
s_2          & 1 & 1 & 1 & 2 & 1 & 1 & 0 & 0 & 2 & 1 
\end{array}
\right],
\nonumber
\\
&
{\bf E}_4 = \left[
\begin{array}{ccccc|cc|cc|cc}
s_3-2s_1-4 & 1 & 1 & 1 & 2 & 0 & 0 & 1 & 1 & 1 & 2  \\
s_2-2s_1-6 & 1 & 1 & 1 & 2 & 1 & 1 & 0 & 0 & 2 & 1 
\end{array}
\right],
\label{C4} \\
&
{\bf E}_5 = \left[
\begin{array}{ccccc|cc|cc|cc|}
s_3-s_1-1 & 1 & 1 & 1 & 1 & 2 & 2 & 1 & 1 & 0 & 0  \\
s_2-s_1-2 & 1 & 1 & 1 & 1 & 1 & 1 & 0 & 0 & 1 & 1 
\end{array}
\right],
\nonumber
\\
&
{\bf E}_6 = \left[
\begin{array}{ccccc|cc|cc|cc|}
s_3-2s_1-4 & 1 & 1 & 1 & 1 & 2 & 2 & 1 & 1 & 0 & 0  \\
s_2-s_1-2 & 1 & 1 & 1 & 1 & 1 & 1 & 0 & 0 & 1 & 1 
\end{array}
\right],
\nonumber
\eea
where the groups of C-columns are bounded from the right by a vertical line.
Recalling that $W({\bf E}_k) \equiv W({\bf s}_k,{\bf D}_k)$
we find three pairs of matrices with coinciding numerical parts ${\bf D}_k$ and different 
argument columns  ${\bf s}_k$ -- namely, ${\bf D}_k = {\bf D}_{k+1}$ for $k=1,3,5$.

Use the recipies from Section \ref{Collinear} to rewrite $W({\bf s}_1,{\bf D}_1)$.
We have ${\bf C}_1 = \{1,1\}^T$ with ${\bf u}_1 = \{1,1,1,2\}$; 
${\bf C}_2 = \{2,1\}^T$ with ${\bf u}_2 = \{1,1\}$, and
${\bf C}_3 = \{1,0\}^T$ with ${\bf u}_3 = \{1,1\}$.
Use (\ref{noncol05}) to write 
\be
{\bf S}_1({\bf k}) = {\bf s}_1-\sum_{p=1}^3 k_p {\bf C}_p = 
\{s_3-k_1-2k_2-k_3,s_2-k_1-k_2\}^T, 
\quad
\bar{\bf D}_{1} = \left(
\begin{array}{cc}
 3 & 0  \\
 2 & 1
\end{array}
\right),
\label{C5}
\ee 
and
\be
 \!\!\!W({\bf s}_1,{\bf D}_1)\!= \!\!\!
\sum_{k_1=0}^{K_1} \sum_{k_2=0}^{K_2} \sum_{k_3=0}^{K_3} 
H({\bf S}_1({\bf k}))W({\bf S}_1({\bf k}),\bar{\bf D}_{1})  W(k_1,\{1,1,1,2\}) W(k_2,\{1,1\}) W(k_3,\{1,1\}). 
\label{C6}
\ee 
The upper limits $K_i$ in (\ref{C6}) are determined by the nonnegativeness of the 
elements of the vector ${\bf S}_1({\bf k})$ in (\ref{C5}).
The relation (\ref{C6}) admits further simplification as the scalar partition 
of any integer number $k$ into a set of $m$ units equals the binomial 
$\binom{k+m-1}{m-1}$, providing $W(k_i,\{1,1\})=\binom{k_i+1}{1} = k_i+1$.
The reduction of the double partition $W({\bf S}_1({\bf k}),\bar{\bf D}_{1})$
leads to
\be
W({\bf S}_1({\bf k}),\bar{\bf D}_{1}) = 
W(s_3-k_1-2k_2-k_3,\{3\}) - W(2s_3-3s_2+k_1-k_2-2k_3-3,\{3\}),
\label{C7}
\ee 
where both scalar partitions have the form $W(s,\{n\})$ that evaluates 
to $W(s,\{n\}) = \delta_{s \bmod{n},0}$ for nonnegative $s$. Collecting these results we reduce 
(\ref{C6}) to the following combination of the scalar partitions
\bea
W({\bf s}_1,{\bf D}_1)  &=& 
\sum_{k_1=0}^{K_1} \sum_{k_2=0}^{K_2} \sum_{k_3=0}^{K_3} H({\bf S}_1({\bf k})) (k_2+1)(k_3+1)
W(k_1,\{1,1,1,2\})
\nonumber\\
&\times&
\left[W(s_3-k_1-2k_2-k_3,\{3\}) - W(2s_3-3s_2+k_1-k_2-2k_3-3,\{3\})\right].
\label{C8}
\eea
The term $W({\bf s}_2,{\bf D}_2) = W({\bf s}_2,{\bf D}_1)$ 
produces the relation similar to (\ref{C6}) with the only change ${\bf S}_1({\bf k}) 
\Rightarrow {\bf S}_2({\bf k})$ as follows
\be
{\bf S}_2({\bf k}) = {\bf s}_2-\sum_{p=1}^3 k_p {\bf C}_p = 
\{s_3-3s_1-9-k_1-2k_2-k_3,s_2-2s_1-6-k_1-k_2\}^T, 
\label{C9}
\ee 
leading to
\be
W({\bf S}_2({\bf k}),\bar{\bf D}_{2}) = 
W(s_3-3s_1-k_1-2k_2-k_3,\{3\}) - W(2s_3-3s_2+k_1-k_2-k_3-3,\{3\}).
\label{C10}
\ee 
Turn to evaluation of $W({\bf s}_3,{\bf D}_3)$ that is similar to $W({\bf s}_1,{\bf D}_1)$.
We have ${\bf C}_1 = \{1,1\}^T$ with ${\bf u}_1 = \{1,1,1,2\}$; 
${\bf C}_2 = \{0,1\}^T$ with ${\bf u}_2 = \{1,1\}$, and
${\bf C}_3 = \{1,0\}^T$ with ${\bf u}_3 = \{1,1\}$.
Use (\ref{noncol05}) to write 
\be
{\bf S}_3({\bf k}) = {\bf s}_3-\sum_{p=1}^3 k_p {\bf C}_p = 
\{s_3-s_1-1-k_1-k_3,s_2-k_1-k_2\}^T, 
\label{C11}
\ee 
and
\be
W({\bf s}_3,{\bf D}_3)= 
\sum_{{\bf k}=0}^{\bf K} (k_2+1)(k_3+1)H({\bf S}_3({\bf k}))
W({\bf S}_3({\bf k}),\bar{\bf D}_{3})  W(k_1,\{1,1,1,2\}) ,
\
\bar{\bf D}_{3} = \left(
\begin{array}{cc}
 1 & 2  \\
 2 & 1
\end{array}
\right),
\label{C12}
\ee 
where the double partition  $W({\bf S}_3({\bf k}),\bar{\bf D}_{3})$ reduces to
\bea
W({\bf S}_3({\bf k}),\bar{\bf D}_{3}) &=& 
W(2s_3-2s_1-s_2-2-k_1+k_2-2k_3,\{3\}) 
\nonumber \\
&-& W(s_3-2s_2-s_1+k_1+2k_2-k_3-4,\{3\}).
\label{C13}
\eea
Similarly we obtain 
\bea
W({\bf S}_4({\bf k}),\bar{\bf D}_{4}) &=& 
W(2s_3-2s_1-s_2-2-k_1+k_2-2k_3,\{3\}) 
\nonumber \\
&-& W(s_3-2s_2+2s_1+k_1+2k_2-k_3+5,\{3\}),
\label{C14}
\eea
that replaces $W({\bf S}_3({\bf k}),\bar{\bf D}_{3})$ in (\ref{C12})
to find explicit expression for $W({\bf s}_4,{\bf D}_4)$.
Finally, consider $W({\bf s}_5,{\bf D}_5)$ with the 
matrix ${\bf D}_5$ void of noncollinear columns. This case is served by the relations
(\ref{noncol08},\ref{noncol09}) derived in Section \ref{Collinear}. 
We choose the vectors ${\bf C}_p$ as follows
${\bf C}_1 = \{2,1\}^T$ with ${\bf u}_1 = \{1,1\}$; 
${\bf C}_2 = \{1,0\}^T$ with ${\bf u}_2 = \{1,1\}$;
${\bf C}_3 = \{0,1\}^T$ with ${\bf u}_3 = \{1,1\}$, and 
${\bf C}_4 = \{1,1\}^T$ with ${\bf u}_4 = \{1,1,1,1\}$.
Use the vectors ${\bf s}_5 = \{s_3-s_1-1,s_2-s_1-2\},\ {\bf s}_6 = \{s_3-2s_1-4,s_2-s_1-2\}$ to find
\be
{\bf S}_5({\bf k}') = {\bf s}_5-\sum_{p=1}^{3} k_p {\bf C}_p 
= \{s_3-s_1-1-2k_1-k_2,s_2-s_1-2-k_1-k_3\}
= \{k_4,k_4\},
\label{C15}
\ee
and the similar relation for ${\bf S}_6({\bf k}')$.
The expression (\ref{noncol09}) turns into
\bea
W({\bf s}_5,{\bf D}_5)\!\!\!\!\!&=& \!\!\!\!\!
\sum_{{\bf k} =0}^{{\bf K}} H({\bf S}_5({\bf k}'))\binom{k_4+3}{3}
\delta_{s_3-s_1-1-2k_1-k_2,k_4}\delta_{s_2-s_1-2-k_1-k_3,k_4}
 \prod_{j=1}^3(k_j+1),
\label{C16}
\\ 
W({\bf s}_6,{\bf D}_6)\!\!\!\!\!&=&\!\!\!\!\! 
\sum_{{\bf k} =0}^{{\bf K}} H({\bf S}_6({\bf k}'))\binom{k_4+3}{3}  
\delta_{s_3-2s_1-4-2k_1-k_2,k_4}\delta_{s_2-s_1-2-k_1-k_3,k_4}
 \prod_{j=1}^3(k_j+1).
\label{C17}
\eea

\section{Partition with collinear columns -- matrix extension algorithm}
\renewcommand{\theequation}{\thesection\arabic{equation}}
\setcounter{equation}{0}
\label{CollinearExample2}

Repeat the computation of the triple partition $W({\bf E})$ defined by the matrix 
in (\ref{C1}) using the matrix extension approach  discussed in Section \ref{Extension}.
The first reduction is not changed as the matrix ${\bf E}$ is free of such columns,
and our goal to perform the reduction of ${\bf E}_i$ given in (\ref{C4}).
As the numerical columns coincide for the matrix pairs ${\bf D}_k = {\bf D}_{k+1},\ k=1,3,5$
we have to consider just three different cases.
Here for the demonstration of the method we consider only ${\bf E}_1$ and extend it to
\be
\tilde{\bf E}_1 = \left[
\begin{array}{cccccccccccc}
s_0 & 1 & -\delta_{1} & -\delta_{2} & -\delta_{3} & -\delta_{4} 
& -\delta_{5} & -\delta_{6} & -\delta_{7} & -\delta_{8} & -\delta_{9} & -\delta_{10}  \\
s_3 &0 & 1 & 1 & 1 & 2 & 2 & 2 & 1 & 1 & 3 & 0  \\
s_2 &0 & 1 & 1 & 1 & 2 & 1 & 1 & 0 & 0 & 2 & 1 
\end{array}
\right],
\label{D1}
\ee
As $\tilde{\bf E}_1$ has three LZ-columns we end up with 
8 nonzero double partitions $W(\tilde{\bf E}_1^i)$. 
\be
W(\tilde{\bf E}_1) = \sum_{i=1}^8 a_i W(\tilde{\bf s}^i_1,\tilde{\bf D}^i_1),
\quad
a_1=a_2=a_3=a_4=a_7=-1,\
a_5=a_6=a_8=1.
\label{D2}
\ee

\bea
&&
\tilde{\bf D}_1^1 = 
\left[
\begin{array}{cccccccccc}
 1 & \delta_{1}-\delta_{2} &
\delta_{1}-\delta_{3} & 2\delta_{1}-\delta_{4} & \delta_{1}-\delta_{5} 
 & \delta_{1}-\delta_{6}  & -\delta_{7}   & -\delta_{8} & 2\delta_{1}-\delta_{9}  & \delta_{10}-\delta_{1} \\
0 & 0 & 0 & 0 & 1 & 1 & 1 & 1 & 1 & 1
\end{array}
\right],
\nonumber \\  
&& 
\tilde{\bf s}_1^1 =
\left[
\begin{array}{c}
s_2 \delta_1 + \delta_{1}-\delta_{10} \\
s_3-s_2-1 
\end{array}
\right],
\label{D3} 
\eea

\bea
&&\tilde{\bf D}_1^2 = 
\left[
\begin{array}{cccccccccc}
1 & -\delta_{1}+\delta_{2} & \delta_{2}-\delta_{3} & 2\delta_{2}-\delta_{4} & \delta_{2}-\delta_{5} 
& \delta_{2}-\delta_{6}  & -\delta_{7}   & -\delta_{8} & 2\delta_{2}-\delta_{9}  & \delta_{10}-\delta_{2} \\
0 & 0 & 0 & 0 & 1 & 1 & 1 & 1 & 1 & 1
\end{array}
\right],
\nonumber \\
&&\tilde{\bf s}_1^2 =
\left[
\begin{array}{c}
s_2 \delta_2 + \delta_{2}-\delta_{10} \\
s_3-s_2-1 
\end{array}
\right],
\label{D4} \\
&&\tilde{\bf D}_1^3 = 
\left[
\begin{array}{cccccccccc}
1 & -\delta_{1}+\delta_{3}
& \delta_{3}-\delta_{2} & 2\delta_{3}-\delta_{4} & \delta_{3}-\delta_{5} 
& \delta_{3}-\delta_{6}  & -\delta_{7}   & -\delta_{8} & 2\delta_{3}-\delta_{9}  & \delta_{10}-\delta_{3} \\
0 & 0 & 0 & 0 & 1 & 1 & 1 & 1 & 1 & 1
\end{array}
\right],
\nonumber \\
&&\tilde{\bf s}_1^3 =
\left[
\begin{array}{c}
s_2 \delta_3 + \delta_{3}-\delta_{10} \\
s_3-s_2-1 
\end{array}
\right],
\label{D5} \\
&&\tilde{\bf D}_1^4 = 
\left[
\begin{array}{cccccccccc}
1 & \delta_{4}-2\delta_{1}
& \delta_{4}-2\delta_{2} & \delta_{4}-2\delta_{3} & \delta_{4}-2\delta_{5} 
& \delta_{4}-2\delta_{6}  & -2\delta_{7}   & -2\delta_{8} & 2\delta_{4}-2\delta_{9}  & \delta_{10}-\delta_{4} \\
0 & 0 & 0 & 0 & 1 & 1 & 1 & 1 & 1 & 1
\end{array}
\right],
\nonumber \\
&&\tilde{\bf s}_1^4 =
\left[
\begin{array}{c}
s_2 \delta_4 + \delta_{4}-2 \delta_{10} \\
s_3-s_2-1 
\end{array}
\right],
\label{D6} \\
&&\tilde{\bf D}_1^5 = 
\left[
\begin{array}{cccccccccc}
1 & \delta_{1}  - \delta_{5}  & \delta_{2}  - \delta_{5}  & \delta_{3}  -\delta_{5}  
& \delta_{4}  - 2 \delta_{5}  & \delta_{5}  - \delta_{6}  & -\delta_{7}  
& -\delta_{8}  & \delta_{9} -2 \delta_{5}   &  \delta_{10}-\delta_{5}   \\
0 & 1 & 1 & 1 & 2 & 0 & 1 & 1 & 1 & 2
\end{array}
\right],
\nonumber \\
&&\tilde{\bf s}_1^5 =
\left[
\begin{array}{c}
6 \delta_{5}  + s_2 \delta_{5}-\delta_{1}  - \delta_{2}  - \delta_{3}  - \delta_{4}   - \delta_{9}  - \delta_{10}  \\
s_3 - 2 s_2 -8
\end{array}
\right],
\label{D7} \\
&&\tilde{\bf D}_1^6 = 
\left[
\begin{array}{cccccccccc}
1 & \delta_{1}  - \delta_{6}  & \delta_{2}  - \delta_{6}  & \delta_{3}  -\delta_{6}  
& \delta_{4}  - 2 \delta_{6}  &  \delta_{6}-\delta_{5}   & -\delta_{7}  
& -\delta_{8}  & \delta_{9}-2 \delta_{6}   & \delta_{10}-\delta_{6}  \\
0 & 1 & 1 & 1 & 2 & 0 & 1 & 1 & 1 & 2
\end{array}
\right],
\nonumber \\
&&\tilde{\bf s}_1^6 =
\left[
\begin{array}{c}
6 \delta_{6}  + s_2 \delta_{6}-\delta_{1}  - \delta_{2}  - \delta_{3}  - \delta_{4}   - \delta_{9}  - \delta_{10} \\
s_3-s_2-8 
\end{array}
\right],
\label{D8} \\
&&\tilde{\bf D}_1^7 = 
\left[
\begin{array}{cccccccccc}
2 & 2 \delta_{1} - \delta_{9} & 2 \delta_{2} - \delta_{9} & 
2 \delta_{3} - \delta_{9} & 
2 \delta_{4} - 2 \delta_{9} & \delta_{9}-2 \delta_{5}  & 
 \delta_{9}-2 \delta_{6} & -2 \delta_{7} & -2 \delta_{8} & 2 \delta_{10}-\delta_{9}  \\
0 & 1 & 1 & 1 & 2 & 1 & 1 & 2 & 2 & 3
\end{array}
\right],
\nonumber \\
&&\tilde{\bf s}_1^7 =
\left[
\begin{array}{c}
6 \delta_{9} + s_2 \delta_{9} -2 \delta_{1} - 2 \delta_{2} - 2 \delta_{3} - 2 \delta_{4}- 2 \delta_{10} \\
s_3-s_2-8 
\end{array}
\right],
\label{D9} \\
&&\tilde{\bf D}_1^8 = 
\left[
\begin{array}{cccccccccc}
1 & \delta_{10}-\delta_{1} & \delta_{10}-\delta_{2} & \delta_{10} -\delta_{3}
& 2 \delta_{10}-\delta_{4}  & \delta_{10}-\delta_{5}  & \delta_{10}-\delta_{6}  
& -\delta_{7} & -\delta_{8} & 2 \delta_{10} -\delta_{9} \\
0 & 1 & 1 & 1 & 2 & 2 & 2 & 1 & 1 & 3
\end{array}
\right],
\nonumber \\
&&\tilde{\bf s}_1^8 =
\left[
\begin{array}{c}
s_2 \delta_{10} \\
s_3
\end{array}
\right],
\label{D10} 
\eea
The first four double partitions (\ref{D3}--\ref{D6}) have similar structure 
with four LZ-columns out of ten that makes the number of the related scalar partitions equal to 24.
The same is true for the two pairs (\ref{D7}, \ref{D8}) and (\ref{D9}, \ref{D10}) that give 
16 and 18 respectively, -- so the total minimal number of the scalar partitions reaches 58.
This estimate is true if all matrices $\tilde{\bf D}_1^i$ are free of NP- and 
C-columns, but in this case it is impossible to satisfy.
The reason is that $\tilde{\bf D}_1^7$ contains several NP-columns for any set of $\bm \delta=\{\delta_1,\delta_2,\ldots,\delta_{10}\}$.
Thus our goal is to make all columns in $\tilde{\bf D}_1^i$ noncollinear; 
it is possible for the following set of conditions
\bea
&& \delta_{1}-\delta_{5} \ne \delta_{1}-\delta_{6} \ne 
-\delta_{7} \ne -\delta_{8} \ne 2\delta_{1}-\delta_{9} \ne \delta_{10}-\delta_{1},
\nonumber \\
&& \delta_{2}-\delta_{5} \ne \delta_{2}-\delta_{6} \ne 
-\delta_{7} \ne -\delta_{8} \ne 2\delta_{2}-\delta_{9} \ne \delta_{10}-\delta_{2},
\nonumber \\
&& \delta_{3}-\delta_{5} \ne \delta_{3}-\delta_{6} \ne 
-\delta_{7} \ne -\delta_{8} \ne 2\delta_{3}-\delta_{9} \ne \delta_{10}-\delta_{3},
\label{D11} \\
&& \delta_{4}-2\delta_{5} \ne \delta_{4}-2\delta_{6} \ne 
-2\delta_{7} \ne -2\delta_{8} \ne 2\delta_{4}2-\delta_{9} \ne \delta_{10}-\delta_{4},
\nonumber \\
&& \delta_{1}-\delta_{5} \ne \delta_{2}-\delta_{5} \ne \delta_{3}-\delta_{5} \ne 
\delta_{4}/2-\delta_{5} \ne \delta_{5}-\delta_{6} \ne
-\delta_{7} \ne -\delta_{8} \ne \delta_{9}-2\delta_{5} \ne (\delta_{10}-\delta_{5})/2,
\nonumber \\
&& \delta_{1}-\delta_{6} \ne \delta_{2}-\delta_{6} \ne \delta_{3}-\delta_{6} \ne 
\delta_{4}/2-\delta_{6} \ne \delta_{6}-\delta_{5} \ne
-\delta_{7} \ne -\delta_{8} \ne \delta_{9}-2\delta_{6} \ne (\delta_{10}-\delta_{6})/2,
\nonumber \\
&& 2\delta_{1}-\delta_{9} \ne 2\delta_{2}-\delta_{9} \ne 2\delta_{3}-\delta_{9} \ne 
\delta_{4}-\delta_{9} \ne \delta_{9}-2\delta_{5} \ne \delta_{9}-2\delta_{6} \ne
-\delta_{7} \ne -\delta_{8} \ne (2\delta_{10}-\delta_{9})/3,
\nonumber \\
&& \delta_{10}-\delta_{1} \ne \delta_{10}-\delta_{2} \ne \delta_{10}-\delta_{3} \ne 
\delta_{10}-\delta_{4}/2 \ne (\delta_{10}-\delta_{5})/2 \ne (\delta_{10}-\delta_{6})/2 \ne
-\delta_{7} \ne -\delta_{8} \ne (2\delta_{10}- \delta_{9})/3.
\nonumber
\eea
The direct computation shows that the following set $\bm \delta=\{3,2,1,1,1,2,4,5,1,6\}$ 
satisfies all the inequalities (\ref{D11}) leading to:
\bea
&&
\tilde{\bf E}_1^1 = 
\left[
\begin{array}{ccccccccccc}
3 s_2 -3 & 1 & 1 & 2 & 5 & 2  & 1  & -4   & -5 & 5  & 3 \\
s_3-s_2-1 & 0 & 0 & 0 & 0 & 1 & 1 & 1 & 1 & 1 & 1
\end{array}
\right],
\nonumber \\ 
&&\tilde{\bf E}_1^2 = 
\left[
\begin{array}{ccccccccccc}
2 s_2 -4 & 1 & -1 & 1 & 3 & 1  & 0  & -4   & -5 & 3  & 4 \\
s_3-s_2-1 & 0 & 0 & 0 & 0 & 1 & 1 & 1 & 1 & 1 & 1
\end{array}
\right],
\nonumber \\ 
&&\tilde{\bf E}_1^3 = 
\left[
\begin{array}{ccccccccccc}
s_2 -5 & 1 & -2 & -1 & 1 & 0  & -1  & -4   & -5 & 1  & 5 \\
s_3-s_2-1 & 0 & 0 & 0 & 0 & 1 & 1 & 1 & 1 & 1 & 1
\end{array}
\right],
\nonumber \\ 
&&\tilde{\bf E}_1^4 = 
\left[
\begin{array}{ccccccccccc}
s_2 -11 & 2 & -5 & -3 & -1 & -1  & -3  & -8   & -10 & 0  & 11 \\
s_3-s_2-1 & 0 & 0 & 0 & 0 & 1 & 1 & 1 & 1 & 1 & 1
\end{array}
\right],
\label{D16} \\
&&\tilde{\bf E}_1^5 = 
\left[
\begin{array}{ccccccccccc}
s_2 -8 & 1 & 2 & 1 & 0 & -1  & -1  & -4   & -5 & -1  & 5 \\
s_3-2s_2-8 & 0 & 1 & 1 & 1 & 2 & 0 & 1 & 1 & 1 & 2
\end{array}
\right],
\nonumber \\ 
&&\tilde{\bf E}_1^6 = 
\left[
\begin{array}{ccccccccccc}
2 s_2 -2 & 1 & 1 & 0 & -1 & -3  & 1  & -4   & -5 & -3  & 4 \\
s_3-2s_2-8 & 0 & 1 & 1 & 1 & 2 & 0 & 1 & 1 & 1 & 2
\end{array}
\right],
\nonumber \\ 
&&\tilde{\bf E}_1^7 = 
\left[
\begin{array}{ccccccccccc}
s_2 -20 & 2 & 5 & 3 & 1 & 0  & -1  & -3   & -8 & -10  & 11 \\
2s_3-3s_2-8 & 0 & 1 & 1 & 1 & 2 & 1 & 1 & 2 & 2 & 3
\end{array}
\right],
\nonumber \\ 
&&\tilde{\bf E}_1^8 = 
\left[
\begin{array}{ccccccccccc}
6s_2 & 1 & 3 & 4 & 5 & 11  & 5  & 4   & -4 & -5  & 11 \\
s_3 & 0 & 1 & 1 & 1 & 2 & 2 & 2 & 1 & 1 & 3
\end{array}
\right].
\nonumber 
\eea
It must be underlined that the choice of $\delta_i$ made above is
not unique, for example, both $\bm \delta=\{3,2,1,1,1,2,4,5,1,7\}$ 
and $\bm \delta=\{3,2,1,1,1,2,4,5,1,9\}$ also work for $W(\tilde{\bf E}_1)$.

We have (when it is required) to 
transform the matrices $\tilde{\bf E}_1^i$ into $\bar{\bf E}_1^i$
where all LZ-columns have positive integer first element 
by applying the column sign inversion procedure (\ref{mCayley118}) to these 
LZ-columns as described in Section \ref{ZeroMulti}. 
For example, $W(\tilde{\bf E}_1^4) = -W(\bar{\bf E}_1^4)$, where
$$
\bar{\bf E}_1^4 = 
\left[
\begin{array}{ccccccccccc}
s_2 -20 & 2 & 5 &3 & 1 & -1  & -3  & -8   & -10 & 0  & 11 \\
s_3-s_2-1 & 0 & 0 & 0 & 0 & 1 & 1 & 1 & 1 & 1 & 1
\end{array}
\right].
$$
Finally we perform the reduction of the double partitions $W(\bar{\bf E}_1^i)$
to obtain the desired expression of $W({\bf E}_1)$ into scalar partitions.
The other matrices ${\bf E}_i, \ 2 \le i \le 6$ are processed similarly.

\end{document}